# BOUNDARIES OF HYPERBOLIC GROUPS

ILYA KAPOVICH AND NADIA BENAKLI

ABSTRACT. In this paper we survey the known results about boundaries of word-hyperbolic groups.

## 1. INTRODUCTION

In the last fifteen years Geometric Group Theory has enjoyed fast growth and rapidly increasing influence. Much of this progress has been spurred by remarkable work of M.Gromov [95], [96] who has advanced the theory of word-hyperbolic groups (also referred to as Gromov-hyperbolic or negatively curved groups). The basic idea was to explore the connections between abstract algebraic properties of a group and geometric properties of a space on which this group acts nicely by isometries. This connection turns out to be particularly strong when the space in question has some hyperbolic or negative curvature characteristics. This led M.Gromov [95] as well as J.Cannon [48] to the notions of a Gromov-hyperbolic (or "negatively curved") space, a word-hyperbolic group and to the development of rich, beautiful and powerful theory of word-hyperbolic groups. These ideas have caused something of a revolution in the more general field of Geometric and Combinatorial Group Theory. For example it has turned out that most traditional algorithmic problems (such as the word-problem and the conjugacy problem), while unsolvable for finitely presented groups in general, have fast and transparent solutions for hyperbolic groups. Even more amazingly, a result of Z.Sela [167] states that the isomorphism problem is also solvable for torsion-free word-hyperbolic groups.

Gromov's theory has also found numerous applications to and connections with many branches of Mathematics. For example, word-hyperbolic groups turned out to be useful for pursuing Thurston's Geometrization Program for 3-manifolds (see [120], [55], [175]). Hyperbolic groups also lie at the heart of the theory of automatic groups (see [77], [10], [80]) and the related progress in Computational Group Theory. Moreover much of the focus in Group Theory, classical Hyperbolic Geometry and the study of Kleinian Groups, Low-dimensional Topology, the theory of lattices in Lie groups and, to a certain extent, even Riemannian Geometry has shifted to studying "large-scale" geometric and topological properties of various structures. This led, for example, to the emergence of a new and dynamic field of Quasi-isometric

Rigidity ( see [164], [127] , [78], [81]) and the renewed prominence of $CAT(k)$ (or Cartan-Alexandrov-Toponogov) Geometry (see [5], [45]).

To every Gromov-hyperbolic space $X$ (and in particular to every word-hyperbolic group) one can associate the so-called "space at infinity" or *boundary $\partial X$*. This boundary turns out to be an extremely useful and fundamental tool for studying hyperbolic groups and hyperbolic spaces. The boundary has a rich structure (topological, dynamical, metric, quasiconformal, measure-theoretic and algebraic) with deep and important theory and numerous applications. For example, the work of B.Bowditch [35] related the so-called Jaco-Shalen-Johansonn (or JSJ) decomposition of hyperbolic groups, originally discovered by Z.Sela [168], to the topological properties of the boundary. The hyperbolic boundary turned out to be a rather ubiquitous object arising naturally in almost every possible approach to hyperbolic groups and spaces.

In this paper we aim to give a broad overview of known results, methods and ideas involved in studying and using the hyperbolic boundary. The focus of the paper is more on the boundaries of hyperbolic groups and to a lesser extent on the (unimaginably vast) subject of boundaries of Gromov-hyperbolic spaces. For this reason many important themes, such as $CAT(-1)$ geometry, harmonic analysis on trees and the classical theory of Kleinian groups, are given only cursory treatment.

Because of the space limitations we will ordinarily not provide precise proofs, but rather give an informal non-technical discussion of the main arguments. For the same reason our reference list is far from complete. However, we do start with a careful presentation of the basic definitions and concepts. The reader should bear in mind that this is an extremely active and rapidly developing subject with many new and exciting developments taking place with great speed. Thus our bibliography is certain to be somewhat out-of-date by the time this paper appears in print. Our hope is that this paper can serve as the first stepping stone for people interested in the subject and that it will be accessible and useful to graduate students and non-specialists rather than only to the people working in the field.

The paper is organized as follows. The basic definitions and background information are given in Section 2. In Section 3 we discuss the quasiconformal and quasi-Möbius structures on the boundary. Section 4 deals with basic aspects of the topological dynamic of the action by a word-hyperbolic group on its boundary. In Section 5 we relate hyperbolic groups to convergence group actions. Section 6 is devoted to the results relating cohomological properties of a hyperbolic group and its boundary. Section 7 contains the analysis of the relationship between local connectivity of the boundary and splittings of hyperbolic groups over finite and virtually cyclic subgroups. In Section 8 we give an overview of hyperbolic groups with zero and one-dimensional boundary. Section 9 deals with the work of J.Cannon on groups with boundary homeomorphic to the two-sphere. Section 11 deals with random walks and representing hyperbolic boundary as the Poisson boundary and the Martin boundary. Section 10 outlines the connections with $C^*$-algebras, Gromov's conformal boundary and the Novikov conjecture. The use of boundary considerations for investigating the subgroup structure of hyperbolic groups (particularly limit sets and the Cannon-Thurston map) is discussed in Section 12. Section 13 gives a brief introduction to the notion of a geodesic flow. In Section 14 we talk in more detail about Pansu-Gromov conformal dimension of the boundary. Sullivan-Patterson measures are outlined in Section 15. In Section 16 we discuss symbolic dynamics and viewing a hyperbolic group and its boundary as an abstract dynamical system. Finally various miscellaneous topics are presented in Section 17.

## 2. Definition and basic properties of the boundary

In this section we will give the definition and basic properties of hyperbolic spaces, hyperbolic groups and their boundaries. For careful proofs and a more detailed discussion the reader is referred to [95], [66], [92], [2], [45].

A metric space $(X, d)$ is said to be *geodesic* if any two points $x, y \in X$ can be joined be a *geodesic segment* $[x, y]$ that is a naturally parameterized path from $x$ to $y$ whose length is equal to $d(x, y)$. Informally speaking a hyperbolic space is a geodesic metric space where all geodesic triangles are thin. That is to say, a geodesic triangle looks more or less like a "tripod". More generally a polygon with geodesic sides in a hyperbolic space looks similar to a tree.

**Definition 2.1** (Hyperbolic space). A geodesic metric space $(X, d)$ is called $\delta$-*hyperbolic* (where $\delta \geq 0$ is some real number) if for any triangle with geodesic sides in $X$ each side of the triangle is contained in the $\delta$-neighborhood of the union of two other sides. That is, if $\Delta$ is a triangle with geodesic sides $\alpha$, $\beta$, $\gamma$ in $X$ then for any $p \in \alpha$ there is $q \in \beta \cup \gamma$ with $d(p, q) \leq \delta$.

A geodesic metric space is said to be *hyperbolic* if it is $\delta$-hyperbolic for some $\delta \geq 0$.

Even though it is possible to define hyperbolicity for a metric space that is not necessarily geodesic (see [95]), for the purposes of this paper we will restrict ourselves to hyperbolic geodesic metric spaces only.

**Example 2.2.** Any geodesic metric space of finite diameter $D$ is $D$-hyperbolic. An $\mathbb{R}$-tree (e.g. the real line $\mathbb{R}$ or a simplicial tree) is 0-hyperbolic. The standard hyperbolic plane $\mathbb{H}^2$ is $log(\sqrt{2}+1)$-hyperbolic. If $M$ is a closed Riemannian manifold of negative sectional curvature then the universal covering space of $M$ (with the induced metric) is hyperbolic.

Another basic example of hyperbolic spaces is provided by the so-called $CAT(k)$-spaces, where $k < 0$. We recall their definition:

**Definition 2.3.** Let $(X, d)$ be a geodesic metric space and let $k \leq 0$. We say that $(X, d)$ is a $CAT(k)$-*space* if the following holds.

Suppose $x, y, z \in X$ and let $p \in [x, y]$, $q \in [x, z]$. Let $\mathcal{M}_k^2$ be the model space of constant curvature $k$, that is $\mathcal{M}_k^2$ is the Riemannian manifold homeomorphic to $\mathbb{R}^2$ and with constant sectional curvature $k$. Let $x', y', z' \in \mathcal{M}_k^2$ be such that $d(x, y) = d(x', y')$, $d(x, z) = d(x', z')$, $d(y, z) = d(y', z')$. Let $p' \in [x', y']$ and $q' \in [x', z']$ be such that $d(x, p) = d(x', p')$ and $d(x, q) = d(x', q')$. Then

$$d(p, q) \leq d(p', q').$$

Thus geodesic triangles in a $CAT(k)$-space are no fatter than comparison triangles in the model space $\mathcal{M}_k^2$. Note that for $k = 0$ the model space $\mathcal{M}_0^2$ is just the Euclidean plane and for $k < 0$ the model space $\mathcal{M}_k^2$ is a scaled version of the hyperbolic plane $\mathbb{H}^2$. As we noted before, if $k < 0$ then any $CAT(k)$-space is hyperbolic. Moreover, geodesics between any two points in such a space are unique. One can also define a space *of curvature* $\leq k$ (where $k \leq 0$) as a geodesic metric space where every point has a convex $CAT(k)$-neighborhood. We refer the reader to the book of M.Bridson and A.Haefliger [45] and the book of W.Ballmann [5] for more information.

**Definition 2.4** (Hyperbolic Group). A finitely generated group $G$ is said to be word-hyperbolic if there is a finite generating set $S$ of $G$ such that the Cayley graph $\Gamma(G, S)$ is hyperbolic with respect to the word-metric $d_S$.

It turns out that if $G$ is a word-hyperbolic group then for *any* finite generating set $S$ of $G$ the corresponding Cayley graph is hyperbolic, although the hyperbolicity constant depends on the choice of $S$. Moreover, word-hyperbolic groups are not just finitely generated by also finitely presentable.

**Example 2.5.** Any finite group $K$ is word-hyperbolic. For any integer $n \geq 1$ the free group $F_n$ of rank $n$ is word-hyperbolic. The fundamental group of the sphere with two handles is word-hyperbolic. If $G$ is a group given by a finite presentation satisfying the $C'(1/6)$-small cancellation condition then $G$ is word-hyperbolic (see [134], [92]). It can also be shown that in a certain precise sense a "generic" finitely presented group is word-hyperbolic [148], [58].

**Definition 2.6** (Gromov Product). Let $(X, d)$ be a metric space. For $x, y, z \in X$ we define

$$(y, z)_x := \frac{1}{2}(d(x, y) + d(x, z) - d(y, z)).$$

We will sometimes call $(y, z)_x$ the *Gromov product* of $y$ and $z$ with respect to $x$.

In hyperbolic metric spaces Gromov product measures how long two geodesics travel close together. Namely if $x, y, z$ are three points in a $\delta$-hyperbolic metric space $(X, d)$ then initial segments of length $(y, z)_x$ of any two geodesics $[x, y]$ and $[x, z]$ are $2\delta$-Hausdorff close. Moreover in this case the Gromov product $(y, z)_x$ approximates within $2\delta$ the distance from $x$ to a geodesic $[y, z]$.

We will say that two geodesic rays $\gamma_1 : [0, \infty) \to X$ and $\gamma_2 : [0, \infty) \to X$ are *equivalent* and write $\gamma_1 \sim \gamma_2$ if there is $K > 0$ such that for any $t \geq 0$

$$d(\gamma_1(t), \gamma_2(t)) \leq K.$$

It is easy to see that $\sim$ is indeed an equivalence relation on the set of geodesic rays. Moreover, two geodesic rays $\gamma_1, \gamma_2$ are equivalent if and only if their images have finite Hausdorff distance. The Hausdorff distance is defined as the infimum of all numbers $H > 0$ such that the image of $\gamma_1$ is contained in the $H$-neighborhood of the image of $\gamma_2$ and vice versa.

The boundary is usually defined as the set of equivalence classes of geodesic rays starting at a base-point, equipped with the compact-open topology. That is to say, two rays are "close at infinity" if they stay Hausdorff-close for a long time. We will now make this notion precise.

**Definition 2.7** (Geodesic boundary). Let $(X, d)$ be a $\delta$-hyperbolic metric space and let $x \in X$ be a base-point.

We define the *relative geodesic boundary* of $X$ with respect to the base-point $x$ as

$$\partial_x^g X := \{\, [\gamma] \,|\, \gamma : [0, \infty) \longrightarrow X \text{ is a geodesic ray with } \gamma(0) = x \,\}.$$

We also define the *geodesic boundary* of $X$ as

$$\partial^g X := \{\, [\gamma] \,|\, \gamma : [0, \infty) \longrightarrow X \text{ is a geodesic ray in } X \,\}.$$

**Convention 2.8.** Let $x \in X$ be a base-point in a hyperbolic metric space $(X, d)$. We say that a sequence $(x_n)_{n \geq 1}$ of points in $X$ *converges to infinity* if

$$\lim_{i, j \to \infty} \inf (x_i, x_j)_x = \infty.$$

It is easy to see that this definition does not depend on the choice of a base-point. We shall also say that two sequences $(x_n)$ and $(y_n)$ converging to infinity are *equivalent* and write $(x_n) \sim (y_n)$ if

$$\lim_{i, j \to \infty} \inf (x_i, y_j)_x = \infty.$$

Once again, one can show that $\sim$ is an equivalence relation on the set of sequences converging to infinity and that the definition of the equivalence does not really depend on the choice of a base-point $x \in X$.

**Definition 2.9** (Sequential boundary). Let $(X, d)$ be a $\delta$-hyperbolic metric space.

We define the *boundary* of $X$ as

$$\partial X := \{\, [(x_n)] \,|\, (x_n)_{n \geq 1} \text{ is a sequence converging to infinity in } X \,\}.$$

One can define some canonical maps between different boundaries. Namely, if $\gamma$ is a geodesic ray in $X$ starting with $x$ then $[\gamma]$ belongs to both $\partial_x^g X$ and $\partial^g X$. We can put $i_x([\gamma]) := [\gamma]$ which defines a map $i_x : \partial_x^g X \to \partial^g X$. Also, if $\gamma$ is a geodesic ray in $X$ then the sequence $\gamma(n)$ converges to infinity. We put $i([\gamma]) := [(\gamma(n))_{n \geq 1}]$ which gives a well-defined map $i : \partial^g X \to \partial X$. Recall that a metric space $(X, d)$ is called *proper* if all closed metric balls in $X$ are compact. Standard compactness considerations imply that for proper hyperbolic spaces all definitions of the boundary coincide.

**Proposition 2.10.** *Let $(X, d)$ be a proper $\delta$-hyperbolic metric space. Then*

1. *For any $x \in X$ the map $i_x : \partial_x^g X \to \partial^g X$ is a bijection.*
2. *The map $i : \partial^g X \to \partial X$ is a bijection.*
3. *For any two non-equivalent geodesic rays $\gamma_1, \gamma_2 : [0, \infty) \to X$ there exists a biinfinite geodesic $\gamma : (-\infty, \infty) \to X$ such that $\gamma|_{[0,\infty)}$ is equivalent to $\gamma_1$ and $\gamma|_{(-\infty,0]}$, after being re-parameterized from $\gamma(0)$, equivalent to $\gamma_2$.*

**Convention 2.11.** Let $(X, d)$ be a $\delta$-hyperbolic metric space. We will say that a geodesic ray $\gamma : [0, \infty) \to X$ *connects* the point $x = \gamma(0) \in X$ to a point $p \in \partial X$ if $p = [(\gamma(n))_{n \geq 1}]$. Similarly, we will say that a biinfinite geodesic $\gamma : (-\infty, \infty) \to X$ *connects* a point $q \in \partial X$ to a point $p \in \partial X$ if $p = [(\gamma(n))_{n \geq 1}]$ and $q = [(\gamma(-n))_{n \geq 1}]$

Thus we see that in a proper hyperbolic metric space any point in $X$ and any point in $\partial X$ can be connected by a geodesic ray. Moreover, any two distinct points at infinity can be connected by a biinfinite geodesic in $X$.

Even though the boundaries were defined as sets of equivalence classes, they carry a natural topology.

**Definition 2.12** (Topology on the geodesic boundary)**.** Let $(X, d)$ be a $\delta$-hyperbolic metric space and let $x \in X$ be a base-point. For any $p \in \partial_x^g X$ and $r \geq 0$ we define the set

$$V(p, r) := \{ q \in \partial_x^g X | \text{ for some geodesic rays } \gamma_1, \gamma_2 \text{ starting at } x \text{ and with}$$
$$[\gamma_1] = p, [\gamma_2] = q \text{ we have } \liminf_{t \to \infty} (\gamma_1(t), \gamma_2(t))_x \geq r \}.$$

Thus $V(p, r)$ consists of equivalence classes of geodesic rays starting at $x$ which stay $2\delta$-close to a geodesic ray representing $p$ for approximately the distance $r$.

We now topologize $\partial_x^g X$ by setting the basis of neighborhoods for any $p \in \partial_x^g X$ to be the collection

$$\{ V(p, r) | r \geq 0 \}.$$

Thus two geodesic rays starting at $x$ are "close at infinity" if they stay $2\delta$-close for a long time.

In a similar fashion we can define a topology on the sequential boundary.

**Definition 2.13** (Topology on the sequential boundary). Let $(X, d)$ be a $\delta$-hyperbolic metric space and let $x \in X$ be a base-point. For any $p \in \partial X$ and $r \geq 0$ we define the set

$$U(p, r) := \{q \in \partial X | \text{ for some sequences } (x_n), (y_n) \text{ with}$$
$$[(x_n)] = p, [(y_n)] = q \text{ we have } \liminf_{i,j \to \infty} (x_i, y_j)_x \geq r\}.$$

We now endow $\partial X$ with a topology by setting the basis of neighborhoods for any $p \in \partial X$ to be the collection

$$\{U(p, r) | r \geq 0\}.$$

It is not hard to show that the resulting topology does not depend on the choice of a base-point $x \in X$.

Moreover the set $X \cup \partial X$ carries a natural topology as well. Namely, for a point $y \in X$ we use the same basis of neighborhoods as in $X$ itself. For a point $p \in \partial X$ we set the basis of neighborhoods of $p$ in $X \cup \partial X$ to be the collection $\{U'(p, r) | r \geq 0\}$, where

$$U'(p, r) := U(p, r) \cup \{y \in X | \text{ for some sequence } (x_n), [(x_n)] = p$$
$$\text{we have } \lim_{t \to \infty} \inf (x_i, y)_x \geq r\}.$$

Again, the topology on $X \cup \partial X$ does not depend on the choice of a base-point $x \in X$ and it agrees with the earlier notion of "convergence to infinity". Namely for a sequence $(x_n)$ in $X$ and a point $p \in \partial X$

$$[(x_n)] = p \text{ if and only if } \lim_{n \to \infty} x_n = p \text{ in } X \cup \partial X.$$

Once again, it turns out that for proper hyperbolic metric spaces all of the above topologies coincide. Moreover, $\partial X$ is compact and can be considered as a compactification of $X$.

**Proposition 2.14.** *Let $(X, d)$ be a proper $\delta$-hyperbolic metric space. Then*

1. *For any $x \in X$ and $y \in X$ the map*

$$(i_y)^{-1} \circ i_x : \partial_x^g X \to \partial_y^g X$$

   *is a homeomorphism. Thus $i_x$ and $i_y$ induce the same topology on $\partial^g X$.*
2. *For any $x \in X$ the map $i \circ i_x : \partial_x^g X \to \partial X$ is a homeomorphism.*
3. *The topological spaces $\partial X$ and $X \cup \partial X$ are compact.*

**Remark 2.15.** If $(X, d)$ is a proper $CAT(k)$-space where $k < 0$ then the boundary $\partial X$ has some particularly good properties. For example, for any two distinct points in $X \cup \partial X$ there is a unique geodesic connecting them.

**Remark 2.16.** Unlike proper hyperbolic spaces, an arbitrary hyperbolic metric space $X$ does not necessarily have the property that a point $x \in X$ and a point $p \in \partial X$ can be connected by a geodesic ray or that any two distinct points in the boundary can be connected by a biinfinite geodesic. However, a point $x \in X$ and a point $p \in \partial X$ can always be connected by a $(1, 10\delta)$-quasigeodesic ray (see [92], [66], [2] for definition and properties of quasigeodesics in hyperbolic spaces). Indeed, suppose $p = [(x_n)] \in \partial X$. Then for any integer $k \geq 1$ there is a point $y_k \in X$ such that the segment $[x, y_k]$ is $2\delta$-close to an initial segment of $[x, x_n]$ of length $k$ for $n \to \infty$. Then the path

$$\alpha = [x, y_1] \cup [y_1, y_2] \cup [y_2, y_3] \cup \ldots$$

is a $(1, 10\delta)$-quasigeodesic such that the sequence $(y_n)$ is equivalent to $(x_n)$ so that $[(y_n)] = p$. A similar argument shows that any two distinct points $p, q \in \partial X$ can be connected by a biinfinite $(1, 20\delta)$-quasigeodesic path in $X$.

Thus, given an arbitrary $\delta$-hyperbolic metric space $X$, one can define the "quasigeodesic" boundary as the set of equivalence classes of $(1, 20\delta)$-quasigeodesic rays in $X$ (where two rays are equivalent if the Hausdorff distance between them is finite). In this case the quasigeodesic boundary with respect to a base-point, the full quasigeodesic boundary and the sequential boundary of $X$ always coincide and have the same topology (even if $X$ is not proper). Moreover, any two distinct points in $X \cup \partial X$ can be connected by a $(1, 20\delta)$-quasigeodesic path in $X$.

So far the boundary was only endowed with the structure of a topological space. However, it turns out that the boundary $\partial X$ of a proper hyperbolic space is always metrizable (although the metric is not canonical).

**Definition 2.17** (Visual metric). Let $(X, d)$ be a $\delta$-hyperbolic proper metric space. Let $a > 1$ and let $x_0 \in X$ be a base-point. We will say that a metric $d_a$ on $\partial X$ is a *visual metric* with respect to the base-point $x_0$ and the visual parameter $a$ if there is $C > 0$ such that the following holds:

1. The metric $d_a$ induces the canonical boundary topology on $\partial X$.
2. For any two distinct points $p, q \in \partial X$, for any biinfinite geodesic $\gamma$ connecting $p$ to $q$ in $X$ and any $y \in \gamma$ with $d(x_0, y) = d(x_0, \gamma)$ we have:
$$\frac{1}{C} a^{-d(x_0, y)} \leq d_a(p, q) \leq C a^{-d(x_0, y)}.$$

**Theorem 2.18.** [66], [92], [45]

Let $(X, d)$ be a proper $\delta$-hyperbolic metric space. Then:

1. *There is $a_0 > 1$ such that for any base-point $x_0 \in X$ and any $a \in (1, a_0)$ the boundary $\partial X$ admits a visual metric $d_a$ with respect to $x_0$.*

2. *Suppose $d'$ and $d''$ are visual metrics on $\partial X$ with respect to the same visual parameter $a$ and the base-points $x_0'$ and $x_0''$ accordingly. Then $d'$ and $d''$ are Lipschitz equivalent, that is there is $L > 0$ such that*
$$d'(p, q)/L \leq d''(p, q) \leq L d'(p, q) \text{ for any } p, q \in \partial X.$$

3. *Suppose $d'$ and $d''$ are visual metrics on $\partial X$ with respect to the visual parameters $a'$ and $a''$ and the base-points $x_0'$ and $x_0''$ accordingly. Then $d'$ and $d''$ are Holder-equivalent. Namely, there is $C > 0$ such that*
$$[d'(p, q)]^\alpha / C \leq d''(p, q) \leq C[d'(p, q)]^\alpha \text{ for any } p, q \in \partial X$$

*where* $\alpha = \ln a'' / \ln a'$.

**Remark 2.19.** The metric on the boundary is particularly easy to understand when $(X, d)$ is a tree (and so is 0-hyperbolic). In this case $\partial X$ is the space of ends of $X$. The parameter $a_0$ from the above proposition is $a_0 = \infty$ here and for any base-point $x_0 \in X$ and any $a > 1$ the visual metric $d_a$ can be given by an explicit formula:

$$d_a(p, q) = a^{-d(x_0, y)} \text{ for any } p, q \in \partial X$$

where $[x_0, y] = [x_0, p) \cap [x_0, q)$ so that $y$ is the bifurcation point for the geodesic rays $[x_0, p)$ and $[x_0, q)$.

In the case of $CAT(k)$ spaces (where $k < 0$), the boundary admits a visual metric with visual parameter $a$ for any $a \in (1, e^{\sqrt{-k}}]$.

Recall that a function $f : X \to Y$ from a metric space $(X, d_X)$ to a metric space $(Y, d_Y)$ is called a *quasi-isometry* if there is $C > 0$ such that:

1. For any $y \in X$ there is $x \in X$ such that $d(y, f(x)) \leq C$ (that is the image $f(X)$ is co-bounded in $Y$); and
2. For any $x, x' \in X$ we have

$$\frac{1}{C} d_X(x, x') - C \leq d_Y(f(x), f(x')) \leq C d_X(x, x') + C.$$

In this situation $X$ and $Y$ are said to be *quasi-isometric*. Quasi-isometry is a central notion in Geometric Group Theory since it captures the fact that two spaces have similar "large-scale geometry". Since hyperbolicity of a metric space is a "large-scale" property, it is not surprising that we have the following important:

**Proposition 2.20.** [66], [92] *Let* $X, Y$ *be proper geodesic spaces and suppose* $f : X \to Y$ *is a quasi-isometry.*

*Then*

1. $X$ *is hyperbolic if and only if* $Y$ *is hyperbolic.*
2. *The map* $f$ *extends to a canonical homeomorphism* $\hat{f} : \partial X \to \partial Y$.

**Definition 2.21** (Boundary of a hyperbolic group). Let $G$ be a word-hyperbolic group. Then for some (and therefore for any) finite generating set $S$ of $G$ the Cayley graph $\Gamma(G, S)$ is $\delta$-hyperbolic with respect to the word-metric $d_S$. Clearly, $\Gamma(G, S)$ is a proper geodesic metric space. We define the *boundary* $\partial G$ of $G$ as

$$\partial G := \partial \Gamma(G, S).$$

Since the change of a generating set induces a quasi-isometry of the Cayley graphs, the topological type of $\partial G$ does not depend on the choice of $S$.

**Convention 2.22** (Geometric action). We will say that a group $G$ acts on a geodesic space $X$ *geometrically*, if the action is isometric, cocompact (that is $X/G$ is compact) and properly discontinuous, that is for any compact $K \subseteq X$ the set $\{g \in G | gK \cap K \neq \emptyset\}$ is finite.

It can be shown that if $G$ acts on $X$ geometrically then $X$ is *proper*, that is all closed metric balls in $X$ are compact. Notice also that a finitely generated group acts geometrically on its own Cayley graph.

The following classical result of J.Milnor provides a very useful tool for computing boundaries of hyperbolic groups:

**Proposition 2.23.** [92], [66], [77], [45] *Let $G$ be a group acting geometrically on a geodesic metric space $(X, d)$. Then the group $G$ is finitely generated, the space $X$ is proper and for any finite generating set $S$ of $G$ and any $x \in X$ the orbit map*

$$t_x : G \to X, \quad t_x : g \mapsto gx$$

*is a quasi-isometry between $(G, d_S)$ and $X$ (where $d_S$ is the word-metric on $G$ corresponding to $S$).*

This fact implies the following theorem which in a certain sense lies at the heart of Gromov's theory of word-hyperbolic groups:

**Theorem 2.24.** [95] *Let $G$ be a group.*
*Then $G$ is word-hyperbolic if and only if $G$ admits a geometric action on a proper hyperbolic metric space $(X, d)$. Moreover, in this case $\partial G$ is homeomorphic to $\partial X$.*

This provides us with a valuable tool for computing boundaries of many word-hyperbolic groups.

**Example 2.25.**

1. If $G$ is finite then $\partial G = \emptyset$.
2. If $G$ is infinite cyclic then $\partial G$ is homeomorphic to the set $\{0, 1\}$ with discrete topology.
3. If $n \geq 2$ and $G = F_n$, the free group of rank $n$ then $\partial G$ is homeomorphic to the space of ends of a regular $2n$-valent tree, that is to a Cantor set.
4. Let $S_g$ be a closed oriented surface of genus $g \geq 2$ and let $G = \pi_1(S_g)$. Then $G$ acts geometrically on the hyperbolic plane $\mathbb{H}^2$ (the universal cover of $S_g$) and therefore the boundary $\partial G$ is homeomorphic to the circle $S^1$ (the boundary of $\mathbb{H}^2$).
5. Let $M$ be a closed $n$-dimensional Riemannian manifold of constant negative sectional curvature and let $G = \pi_1(M)$. Then $G$ is word-hyperbolic and $\partial G$ is homeomorphic to the sphere $S^{n-1}$.
6. Let $G$ be a group acting isometrically and discretely on the standard hyperbolic space $\mathbb{H}^n$ so that the action is convex-cocompact (i.e. there is a convex $G$-invariant subset $Y$ of $\mathbb{H}^n$ such that $Y/G$ is compact). Then $G$ is word-hyperbolic and $\partial G$ is homeomorphic to the limit set $\Lambda G$ of $G$ in the standard compactification $S^{n-1}$ of $\mathbb{H}^n$.
7. Let $H$ be a subgroup of finite index in a group $G$. Then $H$ is word-hyperbolic if and only if $G$ is word-hyperbolic and in that case $\partial H = \partial G$.

The boundary of a hyperbolic group can also be defined as an inverse limit of the sequence of "sphere projections". Namely, suppose $\Gamma(G, S)$ is the Cayley graph of a word-hyperbolic group $G$ with respect to a finite generating set $S$. We fix an order on $S$ and the induced "short-lex" order on the set of geodesic words in $S$. That is if $w$ and $u$ are $S$-geodesic words, we set $w < u$ if either $|w| < |u|$ or if $|w| = |u|$ and $w$ is lexicographically smaller than $u$. For every nontrivial element $g \in G$ let $w_g$ be the smallest, with respect to this order, geodesic representative word for $g$. Let $w'_g$ be the initial segment of $w_g$ of length $|w_g| - 1 = |g|_S - 1$. Let $u_g \in G$ be the element of $G$ represented by $w'_g$. For each integer $n \geq 1$ we define a projection map $\pi_n : S_n \to S_{n-1}$ from the sphere $S_n$ or radius $n$ around 1 in $G$ to the sphere $S_{n-1}$ in $G$ of radius $n - 1$ centered at 1:

$$\pi_n : g \mapsto u_g, \quad g \in S_n.$$

**Theorem 2.26.** [92], [66] *Let $\Gamma(G, S)$ be the Cayley graph of a word-hyperbolic group $G$ with respect to a finite generating set $S$ and let the family of projections $\pi_n : S_n \to S_{n-1}$ be defined as above, where $n \geq 1$.*

*Then $\partial G$ is homeomorphic to the projective limit of the family $(\pi_n)_{n \geq 1}$.*

A similar statement is true for hyperbolic metric spaces with unique geodesics (e.g. $CAT(k)$-spaces where $k < 0$). Theorem 2.26 can be used to construct a locally finite tree $T$ and a continuous finite-to-one map $\pi : \partial T \to \partial G$. Namely, let $G$ be the vertex set of $T$. For $n \geq 1$ and $g \in S_n$ we connect $g$ and $\pi_n(g)$ by an edge in $T$. Then $T$ is a rooted tree with $1 \in G$ as the root. The tree $T$ can be thought of as a subgraph of $\Gamma(G, S)$ and the degree of each vertex in $T$ is bounded by $2\#(S)$. Any end of $T$ corresponds to a geodesic ray in the Cayley graph $\Gamma(G, S)$. This correspondence defines the map $\pi : \partial T \to \partial G$. The above observations can be used to obtain the following:

**Theorem 2.27.** [92], [66] *Let $G$ be an infinite word-hyperbolic group. Then:*

1. *The boundary $\partial G$ has finite topological dimension $\dim(\partial G)$.*
2. *For any visual metric $d_a$ on $\partial G$ the compact metric space $(\partial G, d_a)$ has finite Hausdorff dimension (which is, of course, greater than or equal to $\dim(\partial G)$) .*

One can show using Gromov's "hyperbolic cone" construction that any compact metrizable space can be realized as the boundary of a proper hyperbolic metric space (see [96], [76], [25]). However, the boundaries of hyperbolic groups are much more restricted:

**Theorem 2.28.** *Let $G$ be a word-hyperbolic group. Then exactly one of the following occurs:*

1. *$G$ is finite and $\partial G$ is empty.*
2. *$G$ contains infinite cyclic group as a subgroup of finite index and $\partial G$ consists of two points.*

3. *G contains a subgroup isomorphic to $F_2$ and the boundary $\partial G$ is an infinite perfect (i.e. without isolated points) compact metrizable space.*

If (1) or (2) occurs, $G$ is termed *elementary*. If (3) occurs, $G$ is called *non-elementary*.

## 3. Quasiconformal structure on the boundary and quasi-isometries

We have observed before that even though $\partial G$ is metrizable, the metric on $\partial G$ is not canonical. However there are weaker *quasiconformal* and *quasi-Möbius* structures on $\partial G$ which are canonical.

An *r-annulus* $A$ in a metric space $(X, d)$ (where $r \geq 1$) is the collection of points between two concentric spheres with the ratio of radii equal to $r$:

$$A = \{x \in X | t \leq d(x, x_0) \leq rt\}$$

where $x_0 \in X$ is the center of the annulus and $t > 0$.

The following notion of a quasiconformal map, suggested by P.Pansu [150], reflects the idea that a quasiconformal map should not distort spheres and annuli two much.

**Definition 3.1.** Let $(X, d_X)$ and $(Y, d_Y)$ be metric spaces. We will say that a map $f : X \longrightarrow Y$ is *quasiconformal* if there is a function $\omega : [1, \infty) \longrightarrow [1, \infty)$ with the following property.

For any $r \geq 1$ and any $r$-annulus $A$ in $X$ there is an $\omega(r)$-annulus $A'$ in $Y$ such that $f(A) \subseteq A'$.

Similarly to quasi-conformal structure, one can define the notion of a real cross-ratio and a quasi-Möbius structure on the boundary. Namely, if $a, b, c, d \in \partial G$, F.Paulin [157] defines the *cross ratio* $[a, b, c, d]$ as

$$[a, b, c, d] = \lim \sup_{\substack{a_i \to a, b_i \to b \\ c_i \to c, d_i \to d}} \frac{1}{2}(d(a_i, b_i) + d(c_i, d_i) - d(a_i, d_i) - d(b_i, c_i))$$

(This definition for $CAT(-1)$-spaces was given by J.-P. Otal [149]. In the $CAT(-1)$-case lim sup can be replaced by lim.) Note that cross-ratio is defined with respect to a particular choice of a finite generating set of $G$. Since both ideal and ordinary geodesic quadrilaterals in a hyperbolic space are thin, the above definition is best illustrated by drawing a picture of a finite tree with four vertices $a_i, b_i, c_i, d_i$ and computing the expression under the limit sign in that case. For two spaces $X, Y$ with the cross-ratio structure a map $f : X \to Y$ is said to be *quasi-Möbius* if there is a quasi-Möbius modulus function $\psi : [0, \infty) \to [0, \infty)$ such that:

$$|[f(a), f(b), f(c), f(d)]| \leq \psi(|[a, b, c, d]|) \text{ for any } a, b, c, d \in X.$$

It turns out that while the visual metric on $\partial G$ is non-unique, the quasi-conformal and quasi-Möbius structure on $\partial G$ are in fact canonical and, more-over, these structures are preserved by quasi-isometries. We will say that a homeomorphism $f : E \to E'$ of two metric spaces is a *quasi-conformal equivalence* (respectively *quasi-Möbius equivalence*) if both $f$ and $f^{-1}$ are quasi-conformal (respectively quasi-Möbius). Similarly, two metrics $d$ and $d'$ on a set $E$ are said to be *quasi-conformally equivalent* (respectively *quasi-Möbius equivalent*) if the identity map $i : E \to E$ is a quasi-conformal (respectively quasi-Möbius) equivalence.

**Theorem 3.2.** *(see for example [92], [28])*
*Let $G$ be a word-hyperbolic group. Then:*

1. *For any two visual metrics $d_1$, $d_2$ on $\partial G$ the identity map $id : (\partial G, d_1) \longrightarrow (\partial G, d_2)$ is a quasiconformal equivalence. Moreover, for any two finite generating sets $S_1$ and $S_2$ of $G$, the identity map $id : \partial G \to \partial G$ is a quasi-Möbius equivalence with respect to the quasi-Möbius structures induced by the word metrics $d_{S_1}$ and $d_{S_2}$ on $G$.*

2. *Let $G_1$ and $G_2$ be word-hyperbolic groups and let $f : G_1 \longrightarrow G_2$ be a quasi-isometry (with respect to some choices of word-metrics). Then the induced map $\hat{f} : \partial G_1 \longrightarrow \partial G_2$ is a quasiconformal and a quasi-Möbius homeomorphism.*

We have already seen in Theorem 2.18 that any two visual metrics on the boundary of a hyperbolic space are not just quasiconformally equivalent, but even Holder equivalent in the obvious sense. It turns out that the boundary has an even finer canonical structure called the *power-quasisymmetric* structure [25]. Moreover, a quasi-isometry between two hyperbolic spaces extends to a bi-Holder and power-quasisymmetric equivalence between the boundaries [25].

A natural question arising in this context is if the converse of Theorem 3.2 is true. That is, must two groups with quasiconformally equivalent or quasi-Möbius equivalent boundaries be quasi-isometric?

An important result of F.Paulin [157] and M.Bourdon [28] answers this question in the affirmative:

**Theorem 3.3.** *Let $G_1$ and $G_2$ be word-hyperbolic groups. Suppose $\hat{h} : \partial G_1 \to \partial G_2$ is a homeomorphism which is a quasi-Möbius equivalence or a quasiconformal equivalence. Then $\hat{h}$ extends to a quasi-isometry $h : G_1 \to G_2$.*

Suppose $\hat{h} : \partial G_1 \longrightarrow G_2$ is a quasi-Möbius (or quasiconformal) homeomorphism. There is an obvious candidate for extending $\hat{h}$ "inside" $G_1$ to produce a quasi-isometry $h : G_1 \longrightarrow G_2$. Namely, recall that ideal geodesic triangles in $G_1 \cup \partial G_1$ are $2\delta$-thin and thus each such triangle with three distinct vertices $p, q, r \in G_1 \cup \partial G_1$ has an approximate center $c(p, q, r) \in G_1$. Similar statement holds for $G_2$. Let $g \in G_1$ be an arbitrary element. We can realize $g$ as approximately equal to $c(p, q, r)$ for some $p, q, r \in \partial G_1$. We now put $h(g) := c(\hat{h}(p), \hat{h}(q), \hat{h}(r)) \in G_2$. The challenge is to show that this map is almost correctly defined, that is that up to a finite distance $h(g)$ does not depend on the choices of $p, q, r \in \partial G_1$ with $g \approx c(p, q, r)$. It is here that one must substantially use the fact that $\hat{h}$ preserves the quasi-Möbius (or the quasiconformal) structure.

F.Paulin [157] states the above theorem for spaces more general than Cayley graphs of hyperbolic groups. Therefore some additional assumptions on the space are present in his statement. Later M.Bonk and O.Schramm [25] showed that those additional assumptions are not needed if instead of quasiconformal maps between the boundaries one considers power-quasisymmetric ones. We recall the definition of a quasisymmetric and a power-quasisymmetric map.

**Definition 3.4.** map $f : (X, d_X) \to (Y, d_Y)$ is called a *quasisymmetry* if there is some increasing homeomorphism $\eta : (0, \infty) \to (0, \infty)$ such that for any three distinct points $a, b, c \in X$ we have

$$\frac{d_Y(f(a), f(c))}{d_Y(f(a), f(b))} \leq \eta\big(\frac{d_X(a, c)}{d_X(a, b)}\big)$$

If the above condition is satisfied with

$$\eta(t) = \eta_{\alpha, \lambda}(t) = \begin{cases} \lambda t^{1/\alpha} & \text{for } 0 < t < 1, \\ \lambda t^{\alpha} & \text{for } t \geq 1, \end{cases}$$

for some $\alpha > 0, \lambda \geq 1$ then $f$ is called a *power-quasisymmetry*.

M.Bonk and O.Schramm used a different idea to approach the problem of extending a map between the boundaries to a quasi-isometry between the spaces. Namely, if $Z$ is a compact (or even just bounded) metric space then one can define a *hyperbolic cone* $Cone(Z)$ which is a hyperbolic metric space with $Z$ as the boundary (see also [96], [76]). Essentially it turns out that if the metric on a reasonable $Z$ is replaced by a power-quasisymmetrically equivalent metric, the quasi-isometry type of the cone $Cone(Z)$ is preserved. Moreover, under some mild assumption on a hyperbolic space $X$ the spaces $X$ and $Cone(\partial X)$ are quasi-isometric. This implies Theorem 3.3 (in a considerably more general form). These and other analytically inspired methods were pushed further by M.Bonk, J.Heinonen and P.Koskela in [23] where a number of remarkable results are obtained, including Gromov-hyperbolic analogues of classical uniformization theorems, Gehring-Hayman theorems and so on. There have been some recent significant advances, particularly the remarkable work of J.Cheeger [61], in transferring many classical analytical notions (such as differentiability, tangent spaces, Sobolev spaces, quasiconformality etc) to the setting of arbitrary metric spaces (see also [189, 190, 103]). These results are certain to produce many interesting and substantial applications to Gromov-hyperbolic spaces and their boundaries.

It is worth noting that the topological type of the boundary generally speaking does not determine the quasi-isometry type of the group. For example, on can find two groups $G_1$ and $G_2$ with the following properties. The group $G_1$ acts geometrically on the 2-dimensional complex hyperbolic space and $G_2$ acts geometrically on the 4-dimensional real hyperbolic space. Then $G_1$ and $G_2$ are word-hyperbolic and their boundaries are homeomorphic to the 3-sphere $S^3$. However, it can be shown [150], [32] that the quasiconformal structures on $S^3$ induced by $G_1$ and $G_2$ are not equivalent and hence $G_1$ and $G_2$ are not quasi-isometric. It is also possible [30] to produce two 2-dimensional hyperbolic buildings whose fundamental groups $H_1$ and $H_2$ have the same topological boundary (namely the Menger curve), but $\partial H_1$ and $\partial H_2$ have different conformal dimensions. Thus (see Section 3) $H_1$ is not quasi-isometric to $H_2$.

We refer the reader to [157], [104], [25], [28], [27], [29], [30], [33], [23] for a more detailed discussion of the quasiconformal structure on the boundary of Gromov-hyperbolic spaces.

As we observed before $CAT(-1)$-spaces provide an important subclass of hyperbolic metric spaces. M.Bourdon [26] showed that in this case one can define not just a quasi-conformal, but a canonical conformal structure on the boundary which is compatible with the quasi-Möbius and quasi-conformal structures. If $X$ is a CAT(-1)-space, M.Bourdon [26] constructs a family of visual metrics $d_x$ (where $x \in X$) on the boundary $\partial X$ which is *conformal* in the following sense:

1. For any isometry $g \in Isom(X)$, any $x \in X$ and any $a, b \in \partial X$ we have $d_{gx}(ga, gb) = d_x(a, b)$.
2. For any $x, y \in X$ the metrics $d_x$ and $d_y$ are in the same *conformal class*. More precisely, for any $a, b \in \partial X$, $a \neq b$, we have

$$\frac{d_y(a, b)}{d_x(a, b)} = exp[\frac{1}{2}(B_a(x, y) + B_b(x, y))],$$

where $B_a(x, y) = h_a(x) - h_a(y)$ and where $h_a$ is the Busemann function corresponding to $a \in \partial X$ (see Section 16 for the definition of Busemann functions).

Moreover, it turns out that any isometry of $X$ acts conformally on $\partial X$ in the following sense.

**Theorem 3.5.** [26] *Let $X$ be a proper $CAT(-1)$-space. Then for any $g \in Isom(X)$, any $x \in X$ and any $a \in \partial X$ the limit $\lim_{b \to a} \frac{d_x(ga, gb)}{d_x(b, a)}$ is finite and nonzero. Moreover, this limit is given by the formula*

$$\lim_{b \to a} \frac{d_x(ga, gb)}{d_x(b, a)} = e^{B_a(x, g^{-1}x)}.$$

M.Bonk and O.Schramm [25] used the hyperbolic cone considerations to obtain a number of other extremely useful results regarding hyperbolic spaces and their boundaries. Their main theorem applied to hyperbolic groups implies the following:

**Theorem 3.6.** *Let $G$ be a word-hyperbolic group and let $S$ be a finite generating set of $G$ then*
*(1) The Cayley-graph $\Gamma(G, S)$ is quasi-isometric to a convex subset of $\mathbb{H}^n$ for some $n \geq 2$.*

*(2) After rescaling the metric on $\Gamma(G, S)$ for some visual metric $d$ on the boundary $\partial G$ the space $(\partial G, d)$ admits a bi-Lipschitz embedding into the standard $(n-1)$-sphere.*

A crucial point in establishing the existence of a quasi-isometric embedding in $\mathbb{H}^n$ is proving that the boundary of a hyperbolic group has finite *Assouad dimension.* (See [25] for the definition of Assouad dimension).

## 4. ELEMENTARY TOPOLOGICAL PROPERTIES OF THE BOUNDARY

As we have seen a quasi-isometry between two hyperbolic spaces extends to a homeomorphism between their boundaries. In particular, every isometry of a hyperbolic space is a quasi-isometry and thus induces a homeomorphism of the boundary (which is bi-Lipschitz with respect to any visual metric). It turns out that any isometry belongs to one of three types: *elliptic, parabolic or loxodromic (hyperbolic)*:

**Proposition 4.1** (Classification of isometries)**.** *Let $(X, d)$ be a proper hyperbolic space and let $\gamma : X \to X$ be an isometry of $X$. Then exactly one of the following occurs:*

1. *For any $x \in X$ the orbit of $x$ under the cyclic group $\langle \gamma \rangle$ is bounded. In this case $\gamma$ is said to be elliptic.*
2. *The homeomorphism $\gamma : \partial X \to \partial X$ has exactly two distinct fixed points $\gamma^+, \gamma^- \in \partial X$. For any $x \in X$ the $\langle \gamma \rangle$-orbit map $\mathbb{Z} \to X$, $n \mapsto \gamma^n x$ is a quasi-isometric embedding and $\lim_{n \to \infty} \gamma^n x = \gamma^+$, $\lim_{n \to \infty} \gamma^{-n} x = \gamma^-$. In this case $\gamma$ is said to be hyperbolic or loxodromic.*
3. *The homeomorphism $\gamma : \partial X \to \partial X$ has exactly one fixed point $\gamma^+ \in \partial X$. For any $x \in X$ $\lim_{n \to \infty} \gamma^n x = \gamma^+$, $\lim_{n \to \infty} \gamma^{-n} x = \gamma^+$. In this case $\gamma$ is said to be parabolic.*

Notice that if $\gamma$ has finite order in the isometry group of $X$ then $\gamma$ has to be elliptic.

A word-hyperbolic group $G$ acts by isometries on its Cayley graph and therefore this action extends to the action of $G$ on $\partial G$ by homeomorphisms. We will summarize the basic topological properties of the action of $G$ on $\partial G$ in the following statement.

**Proposition 4.2.** [66], [2], [92], [95]

*Let $G$ be a word-hyperbolic group. Then*

1. *Any element $g \in G$ of infinite order in $G$ acts as a loxodromic isometry of the Cayley graph of $G$. There are exactly two points in $\partial G$ fixed by $g$: the point $g^+ = \lim_{n \to \infty} g^n$ and the point $g^- = \lim_{n \to \infty} g^{-n}$. The points $g^+$ and $g^-$ are referred to as poles or rational points corresponding to the infinite cyclic subgroup $< g > \le G$.*
2. *The action of $G$ on $\partial G$ is minimal, that is for any $p \in \partial G$ the orbit $Gp$ is dense in $\partial G$.*
3. *The set of rational points $Q(G) := \{g^+ | g \in G \text{ is of infinite order } \}$ is dense in $\partial G$.*

4. *The set of pole-pairs $Q'(G) := \{(g^-, g^+) | g \in G$ is of infinite order $\}$ is dense in $\partial G \times \partial G$.*

The above facts are fairly elementary, but the accurate proof does require some care. For example, perhaps the easiest way to see that $Q(G)$ is dense in $\partial G$ is using the automatic structure on $G$ [42]. Namely, the language $L$ of all geodesic words (with respect to some finite generating set $S$ of $G$) is regular and thus recognized by some finite state automaton $M$ (see for example [77]). An element $p \in \partial G$ is realized by a geodesic ray $[1, p)$ in the Cayley graph of $G$ which can be thought of as an infinite geodesic word $w$. Choose a long initial segment $u$ of $w$. Since $w$ can be "read" in $M$ as a label of an infinite path starting from the start state, we can find an $M$-cycle in $w$ which occurs after $u$. That is to say $w$ has an initial segment $uu'v$ such that the word $uu'v^n$ is geodesic for any $n \geq 1$. Then it is easy to see that for $g = (uu')v(uu')^{-1}$ the points $g^+$ and $p$ are close in $\partial G$.

It turns out that the action by an element of infinite order on the boundary of a hyperbolic group has very simple "North-South" dynamics:

**Theorem 4.3.** [66], [2], [92], [95]

*Let $G$ be a word-hyperbolic group and let $g \in G$ be an element of infinite order. Then for any open sets $U, V \subseteq \partial G$ with $g^+ \in U, g^- \in V$ there is $n > 1$ such that $g^n(\partial G - V) \subseteq U$.*

This elementary observation has immediate strong consequences for the topological type of $\partial G$. Namely, it can be used to obtain the following folklore result stating that the boundary of a hyperbolic group is usually a highly non-Euclidean "fractal" space:

**Theorem 4.4.** *Let $G$ be an infinite hyperbolic group. Suppose $\partial G$ contains an open subset homeomomorphic to $\mathbb{R}^n$, where $n \geq 2$. Then $\partial G$ is homeomorphic to $S^n$.*

*Proof.* Since pole-pairs $Q'(G)$ are dense in $\partial G$, there is $g \in G$ of infinite order such that $g^+, g^-$ are both contained in the open subset of $\partial G$ homeomorphic to $\mathbb{R}^n$. Let $U', V'$ with $g^+ \in U', g^- \in V'$ be disjoint open neighborhoods of $g^+$ and $g^-$ which are both homeomorphic to $\mathbb{R}^n$. Thus there are homeomorphisms $f^+ : \mathbb{R}^n \longrightarrow U'$ and $f^- : \mathbb{R}^n \longrightarrow V'$ such that $f^+(0) = g^+$ and $f^-(0) = g^-$. Let $B_2$ and $B_1$ be open balls of radius 2 and 1 centered at the origin in $\mathbb{R}^n$. Put $U = f^+(B_2)$ and $V = f^-(B_1)$. Thus $U$ is an open neighborhood of $g^+$ and $V$ is an open neighborhood of $g^-$ in $\partial G$. Therefore there is $m > 1$ such that $g^m(\partial G - V) \subseteq U$. Consider the sphere $S$ or radius 2 in $\mathbb{R}^n$ centered at the origin. Thus $f^-(S) \subseteq G - V$ is an embedded sphere in $\partial G$. Moreover, this is a collared sphere since we clearly can extend $f^-|_S$ to a topological embedding $S \times [-\epsilon, \epsilon] \longrightarrow \partial G - V$. Consider now the sets $S' = g^m[f^-(S)] \subseteq U \subseteq \partial G$ and $S'' = (f^+)^{-1}(S') \subseteq B_1 \subseteq \mathbb{R}^n$. Then $S''$ is a topologically embedded sphere in $\mathbb{R}^n$ which possesses a bi-collar. Therefore by the Generalized Shoenflies Theorem (see part IV, Theorem 19.11 in [44]) $S''$ separates $\mathbb{R}^n$ into two open connected pieces: one bounded and homeomorphic to an open ball and the other unbounded and homeomorphic to the exterior of a closed ball. Let $C$ be the bounded piece. Since $\partial G - V$ is compact (as a closed subset of $\partial G$), the set $g^m(\partial G - V)$ is compact as well. Hence the set $(f^+)^{-1}(g^m(\partial G - V))$ is a bounded subset of $\mathbb{R}^n$. It is clear that $\partial G - V$ is a closed subset of $\partial G$ with the topological boundary of $\partial G - V$ equal to $f^-(S)$. Since $g^m$ and $f^+$ are homeomorphisms, the interior of $(f^+)^{-1}(g^m(\partial G - V))$ is an open connected subset of $\mathbb{R}^n$ with topological boundary $S''$. We have already seen that $(f^+)^{-1}(g^m(\partial G - V))$ is bounded and therefore it is equal to the closure of $C$. Thus $(f^+)^{-1}(g^m(\partial G - V))$ is homeomorphic to a closed $n$-ball. Hence $\partial G - V$ is also homeomorphic to a closed $n$-ball. The closure of $V$ is equal to $f^-(\bar{B}_1)$ and so it is also a closed $n$-ball by construction. Thus $\partial G$ can be obtained by gluing two closed $n$-balls $\bar{V}$ and $\partial G - V$ along the $(n-1)$-sphere $f^-(S)$. Therefore $\partial G$ is homeomorphic to an $n$-sphere as required. $\qquad\square$

If $\alpha : G \to G$ is an automorphism of a finitely generated group $G$, then $\alpha$ is easily seen to define a quasi-isometry of the Cayley graph of $G$. Thus, if $G$ is word-hyperbolic, $\alpha$ extends to a homeomorphism $\hat{\alpha} : \partial G \to \partial G$. It is interesting to study the dynamics of $\hat{\alpha}$, in particularly the periodic points of $\hat{\alpha}$ (see [133], [132]). G.Levitt and M.Lustig showed [133] that for any word hyperbolic group $G$ "most" automorphisms of $G$ have the same type of simple "North-South" loxodromic dynamics on $\partial G$ as does the translation by an element of infinite order from $G$.

## 5. Convergence group action on the boundary

One of the most important facts about hyperbolic groups is the observation that they act on their boundaries as uniform convergence groups. The notion of convergence group was originally suggested by F.Gehring and G.Martin [91], [90] in order to study the abstract dynamical properties of the action of Kleinian groups on their limit sets (the idea probably goes back to "proximal actions" introduced by H.Furstenberg in [87]). Convergence group actions have been extensively studied in their own right and a powerful machinery for dealing with them has been developed.

**Definition 5.1** (Convergence group). Let $X$ be a compact metrizable space and let $G$ be a group acting on $X$ by homeomorphisms. This action is called a *convergence action* if any of the following equivalent conditions is satisfied.

1. The induced action of $G$ on the space of the so-called "distinct triples" $\Sigma_3(X) := \{(a, b, c) \in X^3 | a \neq b, a \neq c, b \neq c\}$ is properly discontinuous. (Note that $\Sigma_3(X)$ is no longer compact.)

2. For every sequence $(g_n)_{n \geq 1}$ of distinct elements in $G$, there exist a subsequence $g_{n_i}$ and $a, b \in X$ such that $g_{n_i}|_{X - \{a\}}$ converges uniformly on compact subsets to the point $b$. That is to say, for any open neighborhood $U$ of $b$ and any compact subset $C \subseteq X - \{a\}$ there is $i_0$ such that $g_{n_i}(C) \subseteq U$ for any $i \geq i_0$.

A convergence action is called *uniform* if the action of $G$ on $\Sigma_3(X)$ is co-compact.

It turns out that the action of a hyperbolic group on its boundary is a uniform convergence action [82], [137], [187], [193], [39].

Moreover, B.Bowditch [37] proved that in a sense, all uniform convergence actions arise in this fashion.

**Theorem 5.2.** [37] *Let $X$ be a compact metrizable space which is perfect (and hence $X$ is uncountable). Let $G$ be a group acting on $X$ by homeomorphisms. Then the following conditions are equivalent:*

1. *$G$ acts on $X$ as a uniform convergence group.*
2. *$G$ is a non-elementary word-hyperbolic group and $X$ is $G$-equivariantly homeomorphic to $\partial G$.*

One can give another, equivalent, characterization of uniform convergence groups. Suppose $G$ acts as a convergence group on $X$. According to B.Bowditch and P.Tukia, a point $x \in X$ is called a *conical limit point* if there exist $y \in X$ (often equal to $x$), $z \in X$ and a sequence of distinct elements $g_n \in G$ such that $g_n x \to y$ and $g_n|_{X-\{x\}}$ converges uniformly on compact sets to $z$.

**Theorem 5.3.** [37] *Let $G$ act on a compact metrizable perfect space $X$ as a convergence group. Then the following are equivalent:*

1. *$G$ acts on $X$ as a uniform convergence group.*
2. *Every $x \in X$ is a conical limit point.*

One of the first significant results obtained using convergence groups method is the characterization of hyperbolic groups with a circle as their boundary. We have already seen that for the fundamental group $G$ of a closed hyperbolic surface $\partial G$ is homeomorphic to $S^1$ since $G$ is quasi-isometric to the hyperbolic plane $\mathbb{H}^2$. It turns out that this essentially characterizes all groups with $S^1$ as the boundary. A group $G$ is said to be *virtually fuchsian* if $G$ has a geometric action on $\mathbb{H}^2$. The term is justified by the fact that a finitely generated group $G$ is virtually fuchsian if and only if $G$ contains as a subgroup of finite index the fundamental group of a closed hyperbolic surface (see for example [197]).

**Theorem 5.4** (Tukia-Gabai-Freden-Casson-Jungreis). [186], [88], [57],[82]
*Let $G$ be a word-hyperbolic group. Then $\partial G$ is homeomorphic to $S^1$ if and only if $G$ is virtually fuchsian.*

It is worth noting that the actual statement proved in the papers cited above says that a group acts on $S^1$ as a uniform convergence group if and only if the group is virtually fuchsian and the action is topologically conjugate to the action of a virtually fuchsian group on the boundary of $\mathbb{H}^2$. Theorem 5.4 was obtained using deep convergence groups method and is far from trivial. It suffices to say that, as observed by D.Gabai in [88], Theorem 5.4 provides an alternative solution to the Nielsen Realization Problem for surface homeomorphisms (see the paper of S.Kerckhoff [123] for the original proof). Theorem 5.4 also can be used to prove the Seifert fibered space conjecture [88], which forms a significant part of the geometrization program of W.Thurston for 3-manifolds (see [120], [125], [184] for more details).

Using convergence action of a word-hyperbolic group on its boundary allows one to prove many algebraic statements about hyperbolic groups without any use of "hyperbolic" geometry (see for example [91], [90], [35], [37], [84] and other sources). In particular, convergence actions play a central role in Bowditch's approach to the JSJ-decomposition and the study of local connectivity of the boundary (see Section 7 below).

## 6. Cohomological properties of the boundary

The so called *Rips complex*, which can be associated to any hyperbolic group, plays the fundamental role in understanding the cohomologies of hyperbolic groups and their boundaries.

**Definition 6.1** (Rips complex)**.** Let $G$ be a finitely generated group generated by a finite set $S$. Let $d_S$ be the word metric on $G$ corresponding to $S$. For any $d \geq 0$ the simplicial complex $P_d(G, S)$, called the *Rips complex*, is defined as follows. The vertices of $P_d(G, S)$ are elements of $G$. A finite collection of distinct elements $g_0, \ldots, g_k \in G$ spans a $k$-simplex if and only if for all $0 \leq i, j \leq k$ $d_S(g_i, g_j) \leq d$.

It is clear that $P_d(G)$ is a finite-dimensional simplicial complex with a natural left action of $G$. This action is free and transitive on 0-simplices and each $k$-simplex (for $k > 0$) has finite $G$-stabilizer. Thus the action of $G$ on $P_d(G)$ is free if $G$ is torsion-free. Moreover, the quotient $P_d(G)/G$ is compact since $S$ is finite.

A crucial fact regarding hyperbolic groups is (see [2], [92], [45] for accurate proofs):

**Proposition 6.2.** *Let $G$ be a group generated by a finite set $S$ and suppose that the Cayley graph $(\Gamma(G, S), d_S)$ is $\delta$-hyperbolic. Let $d \geq 4\delta + 2$ be an integer. Then $P_d(G, S)$ is contractible.*

Thus $P_d(G)$ can be used to make a free $ZG$-resolution and compute cohomologies of $G$. Standard cohomological considerations now immediately imply the following (see [95], [2], [92], [45]):

**Proposition 6.3.** *Let $G$ be a word-hyperbolic group. Then*

1. *If $G$ is torsion-free then $G$ has finite cohomological dimension.*
2. *If $G$ is virtually torsion-free then $G$ has finite virtual cohomological dimension.*
3. *The group $G$ is of type $FP_\infty$.*
4. *There is $n_0$ such that for all $n \geq n_0$ $H^n(G, \mathbb{Q}) = 0$.*
5. *$H_*(G, \mathbb{Q})$ and $H^*(G, \mathbb{Q})$ are finite-dimensional.*
6. *The group $G$ has only finitely many conjugacy classes of elements of finite order.*

Moreover, as was shown by M.Bestvina and G.Mess [19] one can use the Rips complex $P_d(G)$ to "fill in" $\partial G$ in order to be able to basically view $\partial G$ as the "boundary of a manifold". This trick allows them to connect cohomological properties of $G$ and $\partial G$. Let $X$ be a compact space which is an ANR (absolute neighborhood retract). A closed subset $Q \subseteq X$ is called a *Z-set* or *set of infinite deficiency* if for any open subset $U$ of $X$ the inclusion $U - Q \longrightarrow U$ is a homotopy equivalence. It is easy to see that the boundary of a compact manifold is a $Z$-set in this manifold. We will state the crucial observation for the case of a torsion-free hyperbolic group.

**Proposition 6.4.** [19] *Let $G$ be a torsion-free word-hyperbolic group such that the Cayley graph $(\Gamma(G,S), d_S)$ is $\delta$-hyperbolic with respect to some finite generating set $S$ of $G$. Let $d \geq 4\delta + 2$ be an integer. Then*

1. *$\overline{P}_d(G) = P_d(G) \cup \partial G$ is a compact ANR.*
2. *$\partial G$ is a Z-set in $\overline{P}_d(G)$.*

After proving a similar statement for hyperbolic groups with torsion, M.Bestvina and G.Mess obtain the following corollary using purely topological tools.

**Theorem 6.5.** [19] *Let $G$ be a word-hyperbolic group and let $R$ be a ring with 1. Then:*

1. *We have the isomorphism of $RG$-modules $H^i(G, RG) \cong \check{H}^{i-1}(\partial G, R)$ (the last cohomology is Čech, reduced).*
2. *$dim_R(\partial G) = \max\{n | H^n(G, RG) \neq 0\}$ where the $R$-dimension $dim_R X$ of a topological space $X$ is defined as*

   *$dim_R X := \sup\{n | \check{H}^n(X, A; R) \neq 0 \text{ for some closed subset } A \subseteq X\}$.*

3. *$G$ is a Poincaré duality group of dimension $n$ if and only if $\partial G$ has the integral Čech cohomology of the $(n-1)$-sphere.*
4. *$dim(\partial G) = \max\{n | H^n(G, ZG) \neq 0\}$, where $dim(\partial G)$ is the topological dimension of the boundary $\partial G$;*
5. *If $G$ is torsion-free then $dim(\partial G) = cd(G) - 1$, where*

   *$cd(G) := \max\{m | H^m(G, M) \neq 0 \text{ for some } ZG - \text{module } M\}$*

   *is the cohomological dimension of $G$.*
6. *If $G$ is virtually torsion-free then $dim(\partial G) = vcd(G) - 1$, where $vcd(G)$ is the virtual cohomological dimension of $G$.*

The above statement immediately implies that a generic finitely presented group has one-dimensional boundary. Indeed, suppose $G = < S|R >$ is a torsion-free one-ended group given by a $C'(1/6)$-presentation. Thus $G$ is torsion-free and word-hyperbolic (see [134] for details). Then the presentation complex of $G$ is aspherical and hence the cohomological dimension of $G$ is equal to 2. Therefore $dim(\partial G) = 2 - 1 = 1$.

In the same paper M.Bestvina and G.Mess obtain the following useful application of their methods to 3-manifold groups.

**Theorem 6.6.** [19] *Let $M$ be a closed irreducible 3-manifold with infinite word-hyperbolic fundamental group $G$. Then:*

1. *The universal cover $\tilde{M}$ is homeomorphic to $\mathbb{R}^3$.*
2. *The boundary $\partial G$ is homeomorphic to the 2-sphere $S^2$.*
3. *The compactification $\tilde{M} \cup \partial G$ is homeomorphic to the closed 3-ball.*

Note that if $M$ is a closed irreducible Riemannian 3-manifold with (not necessarily constant) negative sectional curvature, the above theorem implies that the fundamental group of $M$ has $S^2$ as the boundary. This fact plays a central role in Cannon's program of attacking Thurston's Hyperbolization Conjecture (see Section 9).

The approach of M.Bestvina and G.Mess [19] was pursued by M.Bestvina [16] and E.Swenson [179] to obtain more delicate results regarding exotic Steenrod cohomologies of the boundaries of hyperbolic groups.

## 7. BOUNDARY AND SMALL SPLITTINGS OF HYPERBOLIC GROUPS.

If $A, B$ are infinite hyperbolic groups then $G = A * B$ is also hyperbolic and $\partial G$ is disconnected. Indeed, suppose $x \in \partial A$ and $y \in \partial B$ and let $p$ be a path from $x$ to $y$ in $\partial G$. Then $p$ must be approximated by paths $p'$ in $G$ from an element of $A$ to an element of $B$ such that $p'$ stays far away from the identity. But since $G = A * B$, any path from $a \in A$ to $b \in B$ in the Cayley graph of $G$ (with respect to the standard generating set) must pass through 1, which yields a contradiction. Essentially the same argument applies to amalgamated products $G = A *_F B$ and HNN-extensions $G = A*_F$ where $F$ is a finite group. It is also easy to see that there is a map from the set of connected components of $\partial G$ onto the space of ends of $G$. This implies the following:

**Proposition 7.1.** [66], [92], [95] *Let $G$ be non-elementary word-hyperbolic. Then:*

1. *The group $G$ splits (as an amalgamated product or an HNN-extension) over a finite group if and only if $\partial G$ is disconnected*
2. *Suppose $G$ is torsion-free. Then $G$ is freely indecomposable if and only if $\partial G$ is connected.*

Thus it is natural to concentrate on hyperbolic groups with connected non-empty boundary (that is to say one-ended hyperbolic groups). Again, suppose $A$ and $B$ are non-elementary hyperbolic groups and suppose $G = A *_C B$ where $C = <c>$ is an infinite cyclic group. Put $c^+ = \lim_{n\to\infty} c^n \in \partial G$ and $c^- = \lim_{n\to\infty} c^{-n} \in \partial G$. It is easy to see that $c^+$ is a local cut point in $\partial G$. Indeed, suppose $x \in \partial A, y \in \partial B$ are points close to $c^+$. Let $p$ be a path from $x$ to $y$ in $\partial G$ which stays close to $c^+$. We claim that $p$ must in fact pass through $c^+$. Indeed, $p$ is approximated by paths $p'$ from $a \in A$ to $b \in B$ in the Cayley graph of $G$. Since $G = A *_C B$, every such $p'$ must pass through some power $c^i$ of $c$. Moreover we must have $i \to \infty$ since $p'$ approximates $p$. Hence $c^i \to c^+$ and so $p$ passes through $c^+$.

Thus we see that $c^+$ is indeed a local cut point in the boundary. A similar argument applies to hyperbolic groups which split essentially over two-ended (i.e. virtually infinite cyclic) subgroups as an amalgamated product or an HNN-extension. (An amalgamated free product over a two ended subgroup is called *essential* if the associated two-ended subgroup has infinite index in both factors. Similarly, an HNN-extension over a two-ended subgroup is *essential* if the associated two-ended subgroup has infinite index in the base group). A remarkable result of B.Bowditch shows that the converse is also true.

**Theorem 7.2.** [35] *Let $G$ be a one-ended word-hyperbolic group such that $\partial G \neq S^1$. Then*

1. *The boundary $\partial G$ is locally connected and has no global cut points.*
2. *The group $G$ essentially splits over a two-ended subgroup if and only if $\partial G$ has a local cut point.*

This is a very deep, important and beautiful result. Statement (1) in full generality is due to G.Swarup [174] who used convergence groups methods and earlier results of B.Bowditch and G.Levitt to show that if $\partial G$ is connected then $\partial G$ has no global cut-points (Later complete proofs of this fact and some generalizations of it were also given by B.Bowditch [38] and E.Swenson [178]). By the results of M.Bestvina and G.Mess [19] this implies that $\partial G$ is locally connected.

B.Bowditch [35] studies the combinatorics of local cut-points in the boundary of a one-ended hyperbolic group. It turns out that one can use this combinatorics (and that of connected components of the boundary with several local cut-points removed from it) to construct a discrete "pre-tree", that is a structure with a tree-like "betweenness" relation. This pre-tree is used to produce an actual simplicial tree which naturally inherits a $G$-action. The action is then shown to have two-ended edge-stabilizers, which implies Theorem 7.2. Moreover, a detailed analysis of the simplicial tree action of $G$ described above can be used to read-off the JSJ-decomposition of $G$. Roughly speaking, the JSJ decomposition is a canonical splitting of $G$ with two-ended edge stabilizers such that any essential splitting of $G$ over a two-ended subgroup can be obtained from the JSJ using some simple operations. Since Bowditch's construction uses only the topological properties of the boundary, his proof not only provides the uniqueness and existence of the JSJ-decomposition for hyperbolic groups (with or without torsion), but also demonstrates that certain basic properties of the JSJ-decomposition are preserved by quasi-isometries of $G$. Thus, for example, the proof implies that two quasi-isometric (e.g. commensurable) hyperbolic groups have the same number of "quadratically hanging" subgroups in their JSJ-decompositions. The original JSJ-decomposition for torsion-free hyperbolic groups was constructed by Z.Sela [159], [168] by very different methods, namely using Rips' theory of group actions on real trees. Z.Sela proved [168] that for a torsion-free one-ended hyperbolic group $G$ the JSJ-decomposition not only captures all essential splittings over two ended subgroups but also directly encodes the outer automorphism group of $G$. (This is no longer the case for hyperbolic groups with torsion [139].) Several other approaches for obtaining JSJ-decomposition for finitely presented groups have been proposed since then [160], [74], [154], [153], [176].

Theorem 7.2 is the product of deep and extensive research of the subject undertaken by B.Bowditch in a remarkable series of papers [34],[36], [35], [37], [38], [40]. Bowditch's proof of Theorem 7.2 also implies the "algebraic annulus" theorem for hyperbolic groups, which generalizes the classical result of J.Stallings [170].

**Theorem 7.3.** [35] *Let $G$ be a one-ended word-hyperbolic group which is not quasi-fuchsian. Let $C$ be a two-ended subgroup of $G$. Then the number of relative ends $e(G, C) > 1$ if and only if $G$ splits essentially over a two-ended group.*

Indeed, if $e(G, C) > 1$ it is not hard to show that the "endpoints" of $C$ are local cut-points in $\partial G$ and hence the Bowditch's result implies that $G$ splits over a two-ended group. Recently M.Dunwoody and E.Swenson [75] as well as P.Scott and G.Sawrup [166] obtained generalizations of the above theorem to finitely presented groups.

Another result which follows from the purely topological nature of Bowditch's construction is the following:

**Theorem 7.4.** *Let $G_1$ and $G_2$ be quasi-isometric non-elementary one-ended word-hyperbolic groups with the boundaries different from the circle. Then $G_1$ splits over a two-ended subgroup if and only if $G_2$ splits over a two-ended subgroup.*

Recently P.Papasoglu [153] generalized this result to finitely presentabe groups using very different methods.

## 8. BOUNDARY IN LOW DIMENSIONS

As we noted earlier, for a word-hyperbolic group $G$ the topological dimension $dim(\partial G)$ is always finite. Moreover, for torsion-free $G$ we have $dim(\partial G) = cd(G) - 1$.

It turns out that it is easy to classify the groups with zero-dimensional boundary.

**Theorem 8.1.** [95], [66], [92] *Let $G$ be a non-elementary hyperbolic group. Then the following are equivalent.*

1. *$dim \partial G = 0$*
2. *$G$ is virtually free*
3. *$\partial G$ is homeomorphic to the Cantor set.*

The only nontrivial implication is to show that (1) implies (2) and (3). Here $dim \partial G = 0$ implies that $\partial G$ is disconnected and hence $G$ splits over a finite group. Applying the same argument to the vertex groups of this splitting and using the accessibility of finitely presentable groups [73] one can show that this process eventually terminates and so $G$ is obtained from finite groups by iterating free products with amalgamations and HNN-extensions over finite groups. Hence $G$ is virtually free (see for example [8]) and therefore has Cantor set as the boundary.

The situation for groups with $dim(\partial G) = 1$ is substantially more complicated. We have already seen in Theorem 5.4 that groups with $S^1$ as the boundary are precisely virtually fuchsian groups. It is worth noting that most virtually fuchsian groups split over two-ended subgroups (which correspond to simple closed curves on the underlying surfaces). For example, suppose

$$G = \langle a_1, b_1, a_2, b_2 | [a_1, b_1][a_2, b_2] = 1 \rangle$$

is the fundamental group of a sphere with two handles, then $G$ splits as an amalgamated free product of two free groups of rank two over an infinite cyclic group:

$$G = F(a_1, b_1) \underset{[a_1,b_1]=[a_2,b_2]^{-1}}{*} F(a_2, b_2).$$

In general we have the following:

**Theorem 8.2.** [197], [123], [121]

*Let $G$ be a virtually fuchsian hyperbolic group. Then either $G$ splits over a two-ended subgroup or $G$ maps with a finite kernel on the Schwartz triangle group.*

Free products with amalgamations and HNN-extensions over two-ended subgroups, applied to hyperbolic groups with one-dimensional boundary, produce new groups with one-dimensional boundary (provided the result is word-hyperbolic). Thus it is natural to concentrate on one-ended groups with one-dimensional boundary which has no local cut points. Many Kleinian groups of this sort have been long known to have the Sierpinski carpet as their boundary (see for example [11]).

As we observed before, torsion-free small cancellation groups are hyperbolic of cohomological dimension two and thus have 1-dimensional boundary. The topological type of this boundary was not initially understood. N.Benakli [11] in her thesis produced first examples of hyperbolic groups with the boundary homeomorphic to the Menger curve. This was later shown to be the "generic" case by C.Champetier [58] who proved that "most" small cancellation groups have Menger curve as the boundary. These results (as well as the subsequent theorem of M.Kapovich and B.Kleiner [121] mentioned below) rely on some universal properties of the Menger curve among 1-dimensional compact spaces. For example, if $M$ is a connected and locally connected compact metrizable space without local cut-points and of topological dimension one and such that no open nonempty subset of $M$ is planar, then $M$ is homeomorphic to the Menger curve. Thus the crucial point of most proofs concerning the Menger curve is to show that for some reason $\partial G$ is locally non-planar (e.g. by proving that some non-planar graph embeds in an arbitrary small neighborhood of an arbitrary point of $\partial G$).

M.Kapovich and B.Kleiner [121] sharpened the result of C.Champetier by proving the following:

**Theorem 8.3.** *Let $G$ be a one-ended word-hyperbolic group such that $\partial G$ has no local cut points and has topological dimension one. Then one of the following holds:*

1. *$\partial G$ is homeomorphic to the Sierpinski carpet.*
2. *$\partial G$ is homeomorphic to the Menger curve.*

In practice the circle and the Sierpinski carpet are often easy to rule out since they are both planar. Thus it is enough to find a non-planar graph embedded in $\partial G$. This in turn can be done by embedding bigger and bigger copies of this graph in the Cayley graph of $G$ away from the identity (see for example [46] for an illustration of this strategy). Moreover, M.Kapovich and B.Kleiner [121] produced another, cohomological in nature, way of ruling out the Sierpinski carpet:

**Theorem 8.4.** [121] *Let $G$ be a torsion-free hyperbolic group with $\partial G$ homeomorphic to the Sierpinski carpet. Then the Euler characteristic of $G$ is negative.*

We illustrate this point with the following:

**Theorem 8.5.** *Let $G$ be a two-generated one-ended torsion-free $C'(1/6)$-small cancellation group with finite abelianization. Then $G$ is word-hyperbolic and $\partial G$ is homeomorphic to the Menger curve.*

*Proof.* Since $cd(G) = 2$, we know that $dim(\partial G) = 1$. It is clear that $G$ is not virtually Fuchsian since $G$ is non-free, without torsion and it is two-generated. Since $G$ has finite abelianization, it does not split as an HNN-extension. It was observed in [118] that since $G$ is torsion-free word-hyperbolic two-generated and one-ended, $G$ does not essentially split over a two ended subgroup as an amalgamated free product. Thus by the result of Bowditch $\partial G$ has no local cut points.

Therefore by Theorem 8.3 $\partial G$ is either the Menger curve or the Sierpinski carpet. The standard presentation complex of $G$ (corresponding to its finite $C'(1/6)$-presentation) has one vertex, two 1-cells and at least one two-cell (corresponding to the defining relator). Since the presentation is $C'(1/6)$ and $G$ is torsion-free, the presentation complex is in fact a $K(G, 1)$-space (see [134] for the background on small cancellation groups). Hence the Euler characteristic of $G$ is non-negative. Therefore by Theorem 8.4 $\partial G$ cannot be the Sierpinski carpet and so $\partial G$ is homeomorphic to the Menger curve.  $\square$

The results of C.Champetier and Kapovich-Kleiner can also be used to justify the intuitive notion that Cayley graphs of "most" groups are not planar.

**Theorem 8.6.** *Let $G$ be a word-hyperbolic group such that $\partial G$ is homeomorphic to either the Menger curve or the Sierpinski carpet. Then for every finite generating set $S$ of $G$ the Cayley graph $\Gamma(G, S)$ is non-planar.*

*Proof.* We will sketch the argument and leave the details to the reader.

Suppose first $\partial G$ is the Menger curve. Then one can find an embedded copy of the non-planar graph $K_{3,3}$ in $\partial G$. By the definition of $\partial G$ this boundary copy of $K_{3,3}$ can be approximated in $\Gamma(G, S)$ by large copies of $K_{3,3}$ which lie far away from $1 \in G$. Thus $\Gamma(G, S)$ is non-planar.

Suppose now $\partial G$ is the Sierpinski carpet. Then we can find in $\partial G$ an embedded copy of $\Delta$ where $\Delta$ is the 1-skeleton of a tetrahedron. That is $\Delta = K_4$, the complete graph on four vertices (and so $\Delta$ is planar). By the same argument as above we can find a large copy $\Delta'$ of $\Delta$ far away from 1 in the Cayley graph $\Gamma(G, S)$. By adding to $\Delta'$ four geodesics connecting 1 with the four vertices of $\Delta'$ (and rounding the result a little bit around the five vertices) we can find a copy of the cone of $\Delta$ in $\Gamma(G, S)$. Since the cone of $\Delta$ is isomorphic to the complete graph on five vertices $K_5$, it is non-planar. Hence $\Gamma(G, S)$ is also non-planar, as required. $\qquad\square$

## 9. Two-dimensional sphere as the boundary

The famous Thurston's Hyperbolization Conjecture [183], [182], [125], [120] aims to show that a closed atoroidal 3-manifold with infinite fundamental group admits a Riemannian metric of constant negative curvature. The conjecture was settled in a seminal work of Thurston [183], [184] (see also [120]) for Haken manifolds, that is essentially under the assumption that the Poincaré Conjecture holds and that the manifold is aspherical and has an incompressible surface. However, many 3-manifolds are not Haken and a different approach is needed to tackle Thurston's conjecture in general. It turns out that a major chunk of the Hyperbolization Conjecture reduces to the following conjecture of J.Cannon [48]:

**Conjecture 9.1.** *A group $G$ acts geometrically on the 3-dimensional hyperbolic space $\mathbb{H}^3$ if and only if this group is word-hyperbolic and its Gromov boundary is homeomorphic to the 2-sphere.*

The "only if" implication is obvious since in that case $G$ is quasi-isometric to $\mathbb{H}^3$ and hence has $S^2$ as the boundary. The hard part is proving the "if" direction. As we have seen, a similar statement is true for 2-dimensional hyperbolic space (see Theorem 5.4) and false in dimensions higher than three (see Section 14). It is worth noting that a positive solution to the above conjecture would indicate that in dimension three the classes of groups of "variable negative curvature" and "constant negative curvature" essentially coincide. Indeed, if $G$ is the fundamental group of a closed 3-manifold with variable negative curvature, we can still show that $\partial G$ is homeomorphic to $S^2$ (see Theorem 6.5). Conjecture 9.1 would imply that $G$ can also be realized as a fundamental group of a closed hyperbolic manifold (i.e. manifold with constant curvature $-1$). Once again, similar statements fail in higher dimensions [101],[144].

Suppose $G$ is a word-hyperbolic group with the topological 2-sphere as the space at infinity. The challenges in establishing the "if" implication is to show that the action of the group G on the sphere can be made uniformly quasiconformal with respect to the standard conformal structure on the 2-sphere [50]. This is a very formidable difficulty. Indeed, the intrinsic metric on the Cayley graph of G induces a (non-canonical) metric structure on the space at infinity. However, this structure is likely to be very different from the standard one. Thus it usually has Hausdorff dimension greater than two and is ordinarily not smooth. Therefore substantial extra work is needed.

It follows from the results of D.Sullivan [172] that if a group $G$ acts uniformly quasiconformally on the standard two-sphere as a uniform convergence group, then this action is conjugate to a Möbius action (which therefore extends to an isometric action on $\mathbb{H}^3$. P.Tukia [185] proved that a word-hyperbolic group $G$ admits a geometric action on $\mathbb{H}^n$ if and only if $\partial G$ is quasiconformally homeomorphic to the standard $(n-1)$-sphere. Either of these facts can be used to obtain the following:

**Theorem 9.2.** [50] *Let $G$ be a finitely generated group with a finite generating set $S$. Then $G$ admits a geometric action on $\mathbb{H}^3$ if and only if the Cayley graph $\Gamma(G,S)$ is quasi-isometric to $\mathbb{H}^3$.*

J.Cannon has proposed an ingenious and powerful general approach to Conjecture 9.1, which has already produced several remarkable results by J.Cannon and his co-authors (see [48], [52], [49], [51], [53], [54], [55]). Suppose $X$ is a hyperbolic metric space (namely the Cayley graph of a hyperbolic group $G$) with $\partial X$ homeomorphic to $S^2$. Then certain intrinsically defined subsets of $X$ cast disk-type shadows at the sphere at infinity. To be more precise, let $\gamma : [0,\infty) \to X$ be a geodesic ray in $X$. Then for any $t \geq 0$ one can consider the "half-space"

$$H(\gamma, t) := \{x \in X | d(x, \gamma([t,\infty))) \leq d(x, \gamma([0,t]))\}$$

It is easy to see that $H(\gamma, t)$ contains $\gamma([t,\infty))$. This half-space $H(\gamma, t)$ cuts out a "disk" centered at $[\gamma]$ in $\partial X$ defined as:

$$D(\gamma, t) := \{[\gamma'] \in \partial X | \gamma' : [0,\infty) \to X \text{ is a geodesic ray with } \gamma'(0) = \gamma(0)$$
$$\text{such that } \lim_{r \to \infty} d(\gamma'(r), X - H(\gamma, t)) = \infty\}.$$

If we fix a base-point $x_0 \in X$, one can consider the families of disks at infinity

$$\mathcal{D}(n) := \{D(\gamma, n) | \gamma : [0,\infty) \to X \text{ is a geodesic ray with } \gamma(0) = x_0\}.$$

One can then try to understand the internal combinatorics of the way these families of disks ("shingles") cover the 2-sphere. J.Cannon [55] formulated a certain combinatorial axiom (called "the conformality axiom") for these shingling sequences which turns out to guarantee the existence of a uniformly quasiconformal $G$-action and thus would imply the statement of Cannon's conjecture for $G$. The strongest result in this direction so far is stated in a paper by J.Cannon and E.Swenson [55].

**Theorem 9.3.** *A group $G$ has a properly discontinuous co-compact isometric action on $\mathbb{H}^3$ if and only if the following conditions hold:*

1. *The group $G$ is a word-hyperbolic group.*
2. *The boundary $\partial G$ is homeomorphic to $S^2$.*
3. *The sequence of coverings of $S^2$*

$$\mathcal{D}(1), \mathcal{D}(2), \ldots, \mathcal{D}(n), \ldots$$

   *is conformal (see [55] for the precise definition)*

Thus Conjecture 9.1 essentially reduces to the 2-dimensional problem of verifying the conformality of the "shingling sequence" $\mathcal{D}(n)$. Here the biggest challenge is estimating the "combinatorial modulus" for a finite covering of a topological annulus. The main result in this direction is the Combinatorial Riemann Mapping Theorem obtained by J.Cannon in [49] and elaborated on by J.Cannon, W.Floyd and W.Parry in [52].

Recently M.Bonk and B.Kleiner [24] obtained a number of very interesting general results regarding hyperbolic groups with the $n$-sphere as the boundary.

**Theorem 9.4.** [24] *Let $G$ be a word hyperbolic group whose boundary has topological dimension $n$. Assume that $\partial G$ is quasi-conformally homeomorphic to an Ahlfors $n$-regular metric space. Then $G$ acts geometrically on $\mathbb{H}^n$ and $\partial G$ is quasi-conformally homeomorphic to the standard $n$-sphere.*

(The definition of Ahlfors-regular spaces is given in Section 14 below). The above theorem can be used to show, for example, that a word-hyperbolic Coxeter group $G$ with $\partial G$ homeomorphic to the two-sphere acts geometrically on $\mathbb{H}^3$.

## 10. Boundary and $C^*$-algebras

Given any countable connected graph $\Gamma$ (e.g. Cayley graph of a group) the Gelfand-Naimark theory provides standard ways to compactify the vertex set $V = V\Gamma$ using $C^*$-algebras [145], [76].

Namely, let $B = B(\Gamma)$ be the $C^*$-algebra of all bounded $\mathbb{C}$-valued functions on $V$. We will denote the function $V \longrightarrow \mathbb{C}$ defined as $v \mapsto 1, v \in V$ by $\mathbf{1}$. Let $B_0(\Gamma)$ be the subalgebra of $B(\Gamma)$ consisting of functions with finite support. Suppose $A$ (where $B_0 \leq A \leq B$) is a closed subalgebra of $B$ containing $\mathbf{1}$. Then the *spectrum* $sp(A)$ consists of all nonzero characters $A \longrightarrow \mathbb{C}$ which send $\mathbf{1}$ to 1. The spectrum $sp(A)$ has a natural compact-open topology which makes $sp(A)$ a compact topological space. Moreover, if $A$ is separable then $sp(A)$ is metrizable. Namely, if $(f_i)_{i \geq 1}$ is a dense countable sequence in $A$, one can define the metric on $sp(A)$ as

$$d(\tau, \sigma) = \sum_{i=1}^{\infty} \frac{1}{2^i} \min\{|\tau(f_i) - \sigma(f_i)|, 1\}$$

for any $\sigma, \tau \in sp(A)$.

Moreover the vertex set $V$ has a canonical embedding as a dense subset of $sp(A)$:

$$\hat{i} : V \longrightarrow sp(A), \hat{i} : v \mapsto \hat{v} \in sp(A)$$

where

$$\hat{v}(f) := f(v) \text{ for any function } f \in A.$$

Thus $sp(A)$ is a compactification of $\Gamma$ and $\partial_A \Gamma := sp(A) - \hat{V}$ can be regarded as a boundary of $\Gamma$. If $A$ is separable then, as we have seen above, $sp(A)$ is metrizable. In this case the boundary $\partial_A \Gamma$ is sometimes called a *corona* of $\Gamma$ [163]. Different choices of $A$ produce different compactifications of $V$. For example, $A = B(\Gamma)$ corresponds to the Stone-Čech compactification and $A = \langle B_0(\Gamma), \mathbf{1} \rangle$ corresponds to the one-point compactification.

A particularly important and interesting choice of $A$ is provided by the so-called *Gromov-Roe* algebra $\Phi(\Gamma)$. The Gromov-Roe algebra consists of functions on $V$ with rapidly decaying gradient:

$\Phi(\Gamma) := \{f \in B|$ for any $k \geq 1$ there is $C_k > 0$ such that for any two adjacent

vertices $v, u \in V$ with $d(v, v_0) = n$ we have $|f(u) - f(v)| \leq \dfrac{C_k}{(n+1)^k}\}$

(Here $v_0$ is a base-vertex in $\Gamma$ and $d$ is the canonical integer-valued distance function on $V = V\Gamma$). For the choice $A = \Phi(\Gamma)$ the corresponding boundary $\partial_A \Gamma = sp(A) - \hat{V}$ is called the *conformal boundary* of $\Gamma$ [161], [96] and denoted $\partial_\Phi \Gamma$. If $\Gamma$ is a locally finite connected graph with bounded vertex degrees, the Gromov-Roe algebra is easily seen to be separable [76] and so $\partial_\Phi \Gamma$ is a corona. Another important corona in this case is the so-called $l_p$-*corona* $\partial_{\mathcal{L}_p} \Gamma$ corresponding to the subalgebra of $\Phi(\Gamma)$ defined as

$$\mathcal{L}_p(\Gamma) := \Phi(\Gamma) \cap L_p(\Gamma),$$

where $L_p(\Gamma)$ is the space of all real-valued functions on $V$ with $p$-summable gradient. The conformal boundary and the $l_p$-corona play a central role in Gromov's study of conformal dimension of finitely generated groups [96] (see also an excellent paper of G.Elek [76]). It is clear that since $\mathcal{L}_p(\Gamma) \subseteq \Phi(\Gamma)$, the restriction of every unital character on $\Phi(\Gamma)$ to $\mathcal{L}_p(\Gamma)$ is a character of $\mathcal{L}_p(\Gamma)$. Thus there is a natural map between the spectra as well as between the corresponding coronas: $i_p : \partial_\Phi \Gamma \to \partial_{\mathcal{L}_p} \Gamma$. It turns out that for large

$p$ this map is a homeomorphism. This leads M.Gromov to the notion of *conformal dimension* of $\Gamma$:

$$dim\,conf\,\Gamma := \inf\{p > 0 | i_p : \partial_\Phi \Gamma \to \partial_{\mathcal{L}_p}\Gamma \text{ is a homeomorphism}\}.$$

When $\Gamma = \Gamma(G, S)$ is the Cayley graph of a finitely generated group $G$ with respect to a finite generating set $S$, this gives us an important quasi-isometry invariant called the *conformal dimension of $G$*:

$$dim\,conf\,G := dim\,conf\,\Gamma(G, S),$$

(this definition does not depend on the choice of $S$).

A valuable result of J.Roe and N.Higson [107] says that for a word-hyperbolic group the conformal boundary of the Cayley graph is homeomorphic to the standard hyperbolic boundary (in fact the same is true for a Gromov-hyperbolic graph).

**Theorem 10.1.** [107] *Let $G$ be a word-hyperbolic group and let $\Gamma$ be the Cayley graph of $G$ with respect to some finite generating set. Then $\partial G$ is homeomorphic to $\partial_\Phi \Gamma$.*

The homeomorphism $h : \partial G \longrightarrow \partial_\Phi \Gamma, p \in \partial G \mapsto h_p \in sp(\Phi(\Gamma))$ is constructed as follows. Let $p \in \partial G$ and let $\gamma$ be a geodesic ray in $\Gamma$ from $1 \in G$ to $p \in \partial G$. Then for any $f \in \Phi(\Gamma)$ the rapidly decaying gradient condition implies that a finite limit $\lim_{n\to\infty} f(\gamma(n))$ exists. Thus we put $h_p(f) := \lim_{n\to\infty} f(\gamma(n))$ for every $f \in \Phi(\Gamma)$.

The above theorem has immediate interesting applications to the celebrated Novikov conjecture coming from the study of exotic cohomologies in $K$-theory. Namely, J.Roe and N.Higson [107] use Theorem 10.1 to establish the coarse version of the Baum-Connes conjecture for word-hyperbolic groups (see also an article [108] by the same authors). For more information on the Novikov conjecture and related issues the reader is referred to [161],[62], [162], [163], [97], [107], [108].

## 11. Random walks and the Poisson boundary

The classical Poisson formula gives an integral representation of all bounded harmonic functions on the hyperbolic plane in terms of their boundary values on the sphere at infinity. It was generalized to many noncompact symmetric spaces in a classical work of H.Furstenberg [85]. Since then similar results have been obtained for non-positively curved manifolds (see for example [6]). In that context the Poisson formula allows one to study the trajectories of the Brownian motion in the manifold. It turns out that an analogous theory can be developed for random walks on discrete groups (see [115] for a detailed exposition). In fact much of this program can also be carried out for general Markov chains and not just random walks (see [112] for details).

Given a discrete countable group $G$ and a probability measure $\mu : G \longrightarrow [0, 1]$ one can talk about *$\mu$-random walks* on $G$. Namely, the trajectory of a random walk is a sequence

$$g_0, g_1 = g_0 x_1, g_2 = g_0 x_1 x_2, \ldots, g_n = g_0 x_1 x_2 \ldots x_n, \ldots$$

of elements of $G$ corresponding to independent $\mu$-distributed sequences of "jump increments" $x_1, x_2, \ldots$ in $G$.

A function $f : G \longrightarrow \mathbb{C}$ is said to be $\mu$-*harmonic* if it satisfies the following mean-value property:

$$\text{for any } g \in G \quad f(g) = \sum_{x \in G} f(gx)\mu(x).$$

Given a pair $(G, \mu)$ as above one can define a measure space $(B, \nu)$ called the Poisson boundary of $(G, \mu)$. The precise definition is too technical to give here and we refer the reader to [86] and [115] for details. We will give a rather informal idea of what the Poisson boundary is. Imagine that $G$ is embedded in a topological space $X$ so that $X$ is the closure of $G$. Thus the left action of $G$ on itself extends to an action of $G$ on $X$ by homeomorphisms. Let $Y = X - G$ so that $G$ acts on $Y$. Suppose further that almost every trajectory of a $\mu$-random walk on $G$ converges in the topology of $X$ to some point $y \in Y$. Then $Y$ is called a $\mu$-*boundary* of $G$. Note that a $\mu$-boundary $Y$ comes equipped with a *harmonic measure* $\nu$ defined as follows: for a Borel subset $K \subseteq Y$ the measure $\nu(K)$ is defined as the probability that a $\mu$-random walk originating at $1 \in G$ converges to a point in $K$. This measure $\nu$ turns out to be $\mu$-*stationary* in the sense that for any Borel $K \subseteq Y$ we have

$$\nu(K) = \sum_{x \in G} \mu(x)\nu(x^{-1}K).$$

We should stress that our our explanation of the notion of $\mu$-boundary is somewhat misleading. Strictly speaking, a $\mu$-boundary is defined only as a measure space without any topological structure (again, see [115]).

For any $\mu$-boundary $Y$ one can define a linear map from the space of bounded measurable functions on $Y$ to the space $H^\infty(G, \mu)$ of bounded $\mu$-harmonic functions on $G$ via the so-called *Poisson integral map*:

$$\phi : L^\infty(Y, \nu) \longrightarrow H^\infty(G, \mu)$$

$$\phi(f)(g) := \int_Y f(y) dg\nu(y)$$

It is not hard to see that $\phi$ is always injective. The *Poisson boundary* of $(G, \mu)$ is a (unique) $\mu$-boundary $(Y, \nu)$ for which the Poisson integral map $\phi$ is onto (and hence provides an isometry between $L^\infty(Y, \nu)$ and $H^\infty(G, \mu)$). In fact, any $\mu$-boundary can be obtained as an equivariant quotient of the Poisson boundary. We must again stress that the Poisson boundary is simply a measure space with no inherent topological structure.

It turns out that for a very wide class of measures $\mu$ on hyperbolic groups the hyperbolic boundary is a $\mu$-boundary:

**Theorem 11.1.** [111],[114] *Let $G$ be a non-elementary word-hyperbolic group and let $\mu$ be a measure on $G$ whose support generates a non-elementary subgroup of $G$. Then the trajectory of almost every $\mu$-random walk converges to some point in $\partial G$. That is to say, $\partial G$ is a $\mu$-boundary.*

Recall that for a measure $\mu$ on a finitely generated group $G$ the *first moment* $|\mu|$ is defined as

$$|\mu| := \sum_{g \in G} |g| \mu(g)$$

where $|g|$ stands for the geodesic length of $g \in G$ with respect to some fixed finite generating set of $G$. The *first logarithmic moment* of $\mu$ is defined in a similar manner:

$$LM(\mu) := \sum_{g \in G} \mu(g) \log^+ |g|.$$

Another important characteristic of $\mu$ in the *entropy*:

$$H(\mu) := -\sum_{g \in G} \mu(g) \log \mu(g).$$

It is not hard to see that if the first moment of $\mu$ is finite then both the entropy and the first logarithmic moment are finite as well. In [111] V.Kaimanovich shows that for most "reasonable" measures on word-hyperbolic groups the hyperbolic boundary with the harmonic measure is in fact the Poisson boundary.

**Theorem 11.2.** [111] *Let $G$ be a non-elementary word-hyperbolic group and let $\mu$ be a non-degenerate measure with $LM(\mu) < \infty$ and $H(\mu) < \infty$. Then the hyperbolic boundary $\partial G$, when equipped with the harmonic measure $\nu$, is the Poisson boundary of $(G, \mu)$.*

A key geometric ingredient of the proof of the above theorem is the fact that for any two distinct points $p, q \in \partial G$ the intersection of the "convex hull" of $\{p, q\}$ with the ball of radius $n$ in $G$ grows sub-exponentially in $n$. For other results related to the Poisson boundary of hyperbolic groups the reader is referred to [122], [111], [112], [113], [114] [47].

One can also study the so-called *Martin boundary* associated to a random walk $P$ on a graph $\Gamma$ (e.g. the Cayley graph of a word-hyperbolic group). Roughly speaking the Martin boundary $Y$ of $\Gamma$ is a topological space which compactifies $\Gamma$ and which in a certain sense "captures" all *positive* harmonic functions on $\Gamma$ (as opposed to bounded harmonic functions for the Poisson boundary). Every positive harmonic function on $\Gamma$ is represented by a measure on the Martin boundary. The Martin boundary, endowed with the representing measure of the constant 1 function on $\Gamma$, is isomorphic (as a measure space) to the Poisson boundary.

The study of Martin boundaries and positive harmonic functions on graphs is a vast and beautiful field (see for example [192], [3], [124], [7], [194], [112]). For random walks corresponding to finitely supported measures on the Cayley graphs of a word-hyperbolic group (or even on a Gromov-hyperbolic connected graph) the Martin boundary is homeomorphic to the hyperbolic boundary (see [3], [194], [23]). Some applications of these results to Schreier coset graphs of quasiconvex subgroups of hyperbolic groups were recently obtained by I.Kapovich [116]. However, it turns out that Martin boundary of the Cayley graph of a hyperbolic group is less closely connected with the algebraic structure of the group than is the Poisson boundary. For example, it is possible for the Martin boundary of a symmetric random walk on a free group to be homeomorphic to a circle rather than to the Cantor set [7].

## 12. Subgroup structure, the boundary and limit sets

Let $\Gamma$ be a group acting by isometries on a $\delta$-hyperbolic metric space $(X, d)$ (e,g, a word-hyperbolic group acting on its Cayley graph). Let $x \in X$ be a base-point. Then one can define the *limit set* $\Lambda\Gamma \subseteq \partial X$ of $\Gamma$ as the union of all limits in $\partial X$ of sequences of elements from an $\Gamma$-orbit of $x$:

$$\Lambda\Gamma := \{p \in \partial X | p = \lim_{n \to \infty} \gamma_n x \text{ for some sequence } \gamma_n \in \Gamma\}$$

One can also define the *conical limit set* $\Lambda_c\Gamma$ as the collection of those $p \in \partial X$ which are approximated by a sequence from the orbit $\Gamma x$ such that this sequence is contained in a bounded neighborhood of some geodesic ray $r$ with $r(\infty) = p$. Then obviously $\Lambda_c\Gamma \subseteq \Lambda\Gamma$ and both are $\Gamma$-invariant. (For a subgroup $H$ of a hyperbolic group $G$ the notion of a conical limit point of $H$ in $\partial G$ coincides with that given in Section 5 in terms of convergence group action of $H$ on $\partial G$.)

It is clear from the definition that $\Lambda\Gamma$ is a closed $\Gamma$-invariant subset of $\partial X$ and that the definition does not depend on the choice of $x \in X$ . More interestingly, $\Lambda\Gamma$ is the *minimal* closed $\Gamma$-invariant subset of $\partial X$.

**Proposition 12.1.** [63], [65] *Let $\Gamma$ be a group acting by isometries on a hyperbolic space $(X, d)$. Then:*

1. *For any nonempty closed $\Gamma$-invariant subset $A \subseteq \partial X$ we have $\Lambda\Gamma \subseteq A$.*
2. *Suppose the action of $\Gamma$ on $X$ is properly discontinuous. Then $\Gamma$ acts properly discontinuously on $\partial X - \Lambda\Gamma$ and on $(X \cup \partial X) - \Lambda\Gamma$.*

These observations allow one to study $\Lambda H$ using convergence groups methods. Another interesting geometric object associated with the limit set is the *convex hull* $C(\Gamma)$ of $\Gamma$ defined as the union of all biinfinite geodesics in $X$ with both endpoints in $\Lambda\Gamma$. It turns out that $C(\Gamma)$ is always $8\delta$-quasiconvex, where $\delta$ is the hyperbolicity constant of $X$, and $C(\Gamma)$ serves as a good substitute for the classical convex hull of the limit set of Kleinian groups or the concept of minimal $\Gamma$-invariant subtree for groups acting on trees (see [95], [64], [117], [181], [35], [37] and other sources).

In particular, $C(\Gamma)$ can be used to study quasiconvex subgroups of hyperbolic groups and their commensurability properties. Recall that a subset $Y$ of a hyperbolic space $X$ is said to be $\epsilon$-*quasiconvex* in $X$ if for any $y_1, y_2 \in Y$ any geodesic $[y_1, y_2]$ is contained in the $\epsilon$-neighborhood of $Y$. A subgroup $H$ of a word-hyperbolic group $G$ is termed *quasiconvex* if for some (and hence for any) Cayley graph $\Gamma(G, S)$ of $G$ the subgroup $H$ is a quasiconvex subset of $\Gamma(G, S)$. In fact, it turns out that all quasiconvex subgroups of hyperbolic groups are finitely generated and themselves word-hyperbolic. Moreover, a finitely generated subgroup $H$ of a word-hyperbolic group $G$ is quasiconvex in $G$ if and only if for any choices of word metrics $d_H$ on $H$ and $d_G$ on $G$ the inclusion $i : (H, d_H) \to (G, d_G)$ is a quasi-isometric embedding. Here is an assortment of statements that one can prove about quasiconvex subgroups using boundary considerations (compare with the classical results of L.Greenberg for geometrically finite Kleinian groups [93]):

**Theorem 12.2.** [117], [181], [92] *Let $G$ be a word-hyperbolic group acting geometrically on a $\delta$-hyperbolic metric space $(X, d)$, so that $G$ is quasi-isometric to $X$ and is word-hyperbolic. Let $H$ be an infinite subgroup of $G$. Then*

1. *The limit set $\Lambda H$ has either exactly two points or infinitely many points. In the first case $H$ is virtually cyclic and in the second case $H$ contains a free subgroup of rank two.*
2. *The group $H$ acts on $\partial X$ and $\Lambda H$ as a discrete convergence group. For the action of $H$ on $\Lambda H$ a point $p \in \Lambda H$ is a conical limit point in the sense of the definition given in Section 5 if and only if $p \in \Lambda_c H$.*
3. *If $H$ is quasiconvex in $G$ then any point in $\Lambda H$ is a conical limit point, that is $\Lambda H = \Lambda_c H$.*
4. *The group $G$ does not contain an infinite torsion subgroup.*
5. *If $H \leq K \leq G$ and $H$ is normal in $K$ then $\Lambda H = \Lambda K$.*
6. *The subgroup $H$ is quasiconvex in $G$ if and only if $C(H)/H$ has finite diameter.*
7. *If $H$ and $K$ are commensurable subgroups of $G$ then $\Lambda H = \Lambda K$.*
8. *If $H$ is quasiconvex in $G$ then $H$ has finite index in its commensurator*

$$Comm_G(H) = \{g \in G | [H : H \cap g^{-1}Hg] < \infty, [g^{-1}Hg : H \cap g^{-1}Hg] < \infty\}.$$

9. *For any infinite quasiconvex subgroups $H, K \leq G$*

$$\Lambda H \cap \Lambda K = \Lambda(H \cap K).$$

Thus we see that a subgroup of a hyperbolic group is either very small or fairly large. For this reason a subgroup $H$ of a word hyperbolic group is said to be *elementary* if $H$ is either finite or virtually cyclic and *non-elementary* if $H$ contains a free group of rank two.

I.Kapovich and R.Weidmann [119] recently obtained (using boundary considerations) some far-reaching generalizations of the result of M.Gromov-T.Delzant [71], [95] regarding conjugacy classes of non-free subgroups of hyperbolic groups:

**Theorem 12.3.** [119] *Let $G$ be a torsion-free hyperbolic group where all $n$-generated subgroups are quasiconvex. Then $G$ contains only finitely many conjugacy classes of $(n+1)$-generated freely indecomposable non-elementary subgroups.*

The main feature of the proof of Theorem 12.3 is a hyperbolic analog of the classical Nielsen method for studying subgroups of free groups. The hyperbolic version of Nielsen methods deploys looking at minimal networks spanned by convex hulls of several subgroups of a group acting on a hyperbolic space.

Recall that a quasi-isometry between hyperbolic spaces extends to a homeomorphism between their boundaries. Moreover, a quasi-isometric embedding of one hyperbolic space into another extends to a topological embedding on the level of the boundaries. Since quasiconvex subgroups of hyperbolic groups are precisely those finitely generated subgroups which are quasi-isometrically embedded in the ambient hyperbolic group, this immediately implies the following:

**Theorem 12.4.** [66], [92] *Let $H$ be an infinite quasiconvex subgroup of a word hyperbolic group $G$ (so that $H$ is also word-hyperbolic). Then the inclusion of $i : H \rightarrow G$ canonically extends to a topological embedding $\hat{i} : \partial H \rightarrow \partial G$. Moreover, in this case $\hat{i}(\partial H) = \Lambda H$.*

This fact motivates the following notion.

**Definition 12.5.** Let $H$ be an infinite word-hyperbolic subgroup of a word-hyperbolic group $G$. If the inclusion $i : H \rightarrow G$ extends to (necessarily unique) continuous map $\hat{i} : \partial H \rightarrow \partial G$, the map $\hat{i}$ is called the *Cannon-Thurston map* for the pair $H \leq G$. Clearly, if $\hat{i}$ exists, we must have $\hat{i}(\partial H) = \Lambda H$.

A priori if $H$ is not quasiconvex in $G$, it is not clear if a sequence $h_n \in H$ converging to some point in $\partial H$ also converges to a unique point in $\partial G$. Amazingly enough this turns out to be the case for a wide class of non-quasiconvex subgroups. The first example of this sort was provided by J.Cannon and W.Thurston [56]. They showed that if $G$ is the fundamental group of a closed hyperbolic 3-manifold fibering over a circle and if $H$ is the group of the fiber surface then the continuous Cannon-Thurston map $\hat{i} : \partial H \rightarrow \partial G$ exists and is in fact finite-to-one. Moreover $H$ is normal in $G$ and hence the image of $\hat{i}$ equals $\Lambda H = \partial G = S^2$. Since $H$ is a surface group with circle as the boundary, J.Cannon and W.Thurston referred to this map $\hat{i}$ as a "space-filling curve". Later on Mahan Mitra generalized this result in the following two directions.

**Theorem 12.6.** [141] *Suppose we have a short exact sequence of infinite word-hyperbolic groups*

$$1 \rightarrow H \rightarrow G \rightarrow Q \rightarrow 1.$$

*Then a continuous Cannon-Thurston map $\hat{i} : \partial H \rightarrow \partial G$ exists.*

The original Cannon-Thurston example fits the above pattern since in their case $G$ is a semi-direct product of a closed surface group $H$ and an infinite cyclic group $\mathbb{Z}$ along an automorphism of $H$. It is worth noting that by considering the JSJ-decomposition of $H$ one can show that if in the above theorem $H$ is torsion-free then it must be isomorphic to a free product of a free group and several surface groups.

**Theorem 12.7.** [142] *Let* **A** *be a finite graph of groups where all vertex and edge groups are word-hyperbolic and all edge groups are quasi-isometric embeddings. Let* $G = \pi_1(\mathbf{A}, v_0)$ *and suppose that* $G$ *is word-hyperbolic. Then for any vertex group* $H$ *of* **A** *the inclusion* $H \to G$ *induces a continuous Cannon-Thurston map* $\hat{i} : \partial H \to \partial G$.

A key step in establishing both these theorems is to construct, given an $H$-geodesic segment $\lambda$, a certain set $B_\lambda$ containing $\lambda$ in the Cayley graph $\Gamma(G)$ and then define a Lipschitz retraction $\pi_\lambda : \Gamma(G) \to B_\lambda$. This fact is then used to show that $B_\lambda$ is quasiconvex in $G$ (with a quasiconvexity constant independent of $\lambda$). Hence $B_\lambda$ can be used to show that if $\lambda$ was far away from 1 in the $H$-metric, it is still far away from 1 in the $G$-metric which implies the existence of a Cannon-Thurston map. In a subsequent paper [140] M.Mitra developed a theory of "ending laminations" in the context of Theorem 12.7 and proved results similar to those known for 3-dimensional Kleinian groups. A rather intriguing and still open question is whether the Cannon-Thurston map exists for any infinite hyperbolic subgroup of a bigger hyperbolic group.

We must note that an alternative proof of the original Cannon-Thurston result was also obtained by Erica Klarreich [126]. She relied on very different considerations from those of M.Mitra.

## 13. Geodesic flow

Suppose $X$ is the universal covering of a closed Riemannian manifold $X_0$ with negative sectional curvature. Then $X$ is a hyperbolic metric space (even $CAT(k)$ for some $k < 0$) with some particularly good properties of the ideal boundary. Namely, if $\gamma_1, \gamma_2 : [0, \infty) \to X$ are two geodesic rays with $\gamma_1(\infty) = \gamma_2(\infty) \in \partial X$, then the rays $\gamma_1, \gamma_2$ are not just Hausdorff close but converge exponentially fast. Since $X$ is $CAT(k)$ for some $k < 0$, any two distinct points $a, b \in \partial X$ can be connected by a biinfinite geodesic $\gamma : (-\infty, \infty) \to X$ which is unique up-to a translation of $\mathbb{R}$.

Therefore, the space of biinfinite geodesic lines in $X$ (considered up to a translation in $\mathbb{R}$) can be canonically identified with

$$\partial^2 X := \{(a, b) \in \partial X \times \partial X \mid a \neq b\}$$

For a general hyperbolic space $X$ geodesics between two points in the boundary need not be unique and asymptotic geodesics do not have to converge fast and may in fact remain a bounded distance apart. Nevertheless, by analogy with the Riemannian manifold case, M.Gromov [95] suggested how to equip every word-hyperbolic group $G$ with a hyperbolic space quasi-isometric to $G$ such that any two points in the boundary of this space can be connected by a unique geodesic. The details of this construction were made precise by Ch.Champetier [60] and F.Matheus [138].

Let $(X, d)$ be a metric space. We denote by $\mathcal{G}X$ the set

$$\{\gamma : \mathbb{R} \to X \mid \gamma \text{ is a geodesic }\}$$

and equip $\mathcal{G}X$ with the metric $d_{\mathcal{G}}$ as follows:

$$d_{\mathcal{G}}(\gamma_1, \gamma_2) = \int_{-\infty}^{\infty} d(\gamma_1(t), \gamma_2(t)) 2^{-|t|} \, dt.$$

Since $d(\gamma_1(t), \gamma_2(t)) \leq 2t + d(\gamma_1(0), \gamma_2(0))$, the above integral converges. It is not hard to see that $d_{\mathcal{G}}$ is indeed a metric on $\mathcal{G}X$.

The group of isometries of $X$ obviously acts on $\mathcal{G}X$ via composition. The group $\mathbb{Z}_2$ acts by isometries on $\mathcal{G}X$ via the involution $\gamma(t) \longmapsto \gamma(-t)$.

The *geodesic flow* on $X$ is the natural isometric $\mathbb{R}$-action by translation on the elements of $\mathcal{G}X$, given by $\gamma(t) \longmapsto \gamma(t + T)$ for $T \in R$. If $X$ is the universal cover of a closed negatively curved manifold, then the map $\mathcal{G}X \to SX$, $\gamma \mapsto \gamma'(0)$ defines a homeomorphism between $\mathcal{G}X$ and the unit tangent bundle $SX$.

Assume now that $X$ is a proper hyperbolic space. We consider the surjective and continuous projection

$$D : \mathcal{G}X \to \partial^2 X$$

which associates to every biinfinite geodesic of $X$ its endpoints in $\partial X$. In general, for $(a, b) \in \partial^2 X$ the pull-back $D^{-1}(a, b)$ consists of more than one $\mathbb{R}$-orbit since geodesics in $X$ are not unique. The idea of M.Gromov is to identify these geodesics to a single one in a convenient way so that the isometry group of $X$ induces an isometric action on the quotient space.

**Theorem 13.1.** *Let $G$ be a hyperbolic group. There exists a proper metric space $\hat{\mathcal{G}}$ called the geodesic flow space  of $G$ satisfying the following properties:*

1. *The groups $G$ and $\mathbb{Z}_2$ act by isometries on $\hat{\mathcal{G}}$ and the two actions commute. For every $\gamma \in \hat{\mathcal{G}}$, the orbit map $G \to \hat{\mathcal{G}}$, $g \mapsto g \cdot \gamma$ is a quasi-isometry which continuously extends to a canonical $G$-equivariant homeomorphism $\partial G \to \partial \hat{\mathcal{G}}$.*

2. *There is a free $\mathbb{R}$-action on $\hat{\mathcal{G}}$, called the geodesic flow on $\hat{\mathcal{G}}$, which commutes with the action of $G$ and anticommutes with the action of $\mathbb{Z}_2$.*

3. *For $t \in \mathbb{R}$ and $\gamma \in \hat{\mathcal{G}}$ we denote the result of $t$ acting on $\gamma$ by $\gamma(t)$. Then every orbit $\mathbb{R} \to \hat{\mathcal{G}}$, given by $t \mapsto \gamma(t)$ for a fixed $\gamma \in \hat{\mathcal{G}}$, is a quasi-isometric embedding and the orbit space $\hat{\mathcal{G}}/\mathbb{R}$ is homeomorphic to $\partial^2 G$. There is a unique orbit of the geodesic flow that connects two distinct points of $\partial G$. Furthermore there is a homeomorphism between $\hat{\mathcal{G}}$ and $\partial^2 G \times \mathbb{R}$ which conjugates the actions of $\mathbb{R}$.*

4. *Given two points of $\hat{\mathcal{G}}$ such that their orbits by the geodesic flow converge to the same point of $\partial G$, then the orbits get closer to each other exponentially. More precisely, if two $\mathbb{R}$-orbits $\gamma_1(t)$ and $\gamma_2(t)$ satisfy $\gamma_1(\infty) = \gamma_2(\infty) \in \partial G$ then the distance from $\gamma_1(t)$ to $\overline{\gamma_2} = \underset{t}{\cup}\, \gamma_2(t) \subset \hat{\mathcal{G}}$ exponentially decays for $t \to \infty$.*

The construction of $\hat{\mathcal{G}}$ is not canonical but the space $\hat{\mathcal{G}}$ is unique in the following sense:

**Theorem 13.2.** *Let $G$ be a word-hyperbolic group. If $\hat{\mathcal{G}}_1$ and $\hat{\mathcal{G}}_2$ are two proper hyperbolic spaces satisfying the above properties (1), (2), (3) and (4) in Theorem 13.1 then there is a $G \times \mathbb{Z}_2$-equivariant quasi-isometry between $\hat{\mathcal{G}}_1$ and $\hat{\mathcal{G}}_2$ such that every $\mathbb{R}$-orbit of $\hat{\mathcal{G}}_2$ is the homeomorphic image of an $\mathbb{R}$-orbit of $\hat{\mathcal{G}}_1$.*

*Finally, if the group $G$ acts geometrically on a proper hyperbolic space $X$ then there exists a continuous quasi-isometric $G \times \mathbb{Z}_2$-equivariant map $\mathcal{G}X \to \hat{\mathcal{G}}$ which homeomorphically maps every $\mathbb{R}$-orbit of $\mathcal{G}X$ onto an $\mathbb{R}$-orbit of $\hat{\mathcal{G}}$.*

Thus the geodesic flow space $\hat{\mathcal{G}}$ of a group $G$ and its geodesic flow are well defined up to $G$-equivariant quasi-isometries that preserve the orbits. When $G$ is the fundamental group of a negatively curved closed Riemannian manifold $M$ then the space $\hat{\mathcal{G}}$ plays the role of the unit tangent bundle $S\tilde{M}$ of the universal cover $\tilde{M}$ of $M$. In fact, in this case, the geodesic flow has stronger uniqueness properties than those stated in the above theorem (see M.Gromov [95, 98] and F.Matheus [138].)

There is a more general situation when the geodesic flow on a hyperbolic group can be computed explicitly and is well-behaved. Namely, assume that a group $G$ acts by isometries quasi-convex cocompactly (see Section 15 for the definition) on a proper $CAT(-1)$ space $(X, d)$ and hence $G$ is word-hyperbolic . Let $\Lambda G$ be the limit set of $G$ in $\partial X$. We recall that there is a canonical $G$-equivariant quasi-conformal homeomorphism $\partial G \to \Lambda G$.

Under the above assumptions we define

$$\mathcal{G}_G X := \{\gamma \in \mathcal{G}X \,|\, \gamma(\infty) \in \Lambda G, \gamma(-\infty) \in \Lambda G\}.$$

The restriction of $d_{\mathcal{G}}$ to $\mathcal{G}_G X$ will still be denoted by $d_{\mathcal{G}}$. Clearly, $\mathcal{G}_G X$ inherits an $\mathbb{R}$-action and a $\mathbb{Z}_2$-action from $\mathcal{G}X$. The group $G$ acts properly discontinuously and by isometries on $(\mathcal{G}_G X, d_{\mathcal{G}})$. Moreover, it is not hard to see that the geodesic flow space $\hat{\mathcal{G}}$ of $G$ can be equivariantly identified with $\mathcal{G}_G X$. It this situation it is often interesting to consider the quotient-space $E_G X := \mathcal{G}_G X / G$. This "quotient flow space" obviously inherits the $\mathbb{R} \times \mathbb{Z}_2$-action. If $X$ is the universal cover of a closed manifold $X_0$ with sectional curvatures $\leq -1$ and $G = \pi_1(X_0)$ then $E_G X$ can be naturally identified with the unit tangent bundle $SX_0$ of $X_0$.

In this case, M.Bourdon in [28] shows that there is a link between the conformal structure of the limit set and the geodesic flow. Namely, suppose we are given two quasi-convex cocompact actions by isometries of $G$ on two proper $CAT(-1)$ spaces $(X_1, d_1)$ and $(X_2, d_2)$. Let $\Lambda_i G \subseteq \partial X_i$ and $\mathcal{G}_G X_i$ be the corresponding limit sets and geodesic flow spaces, $i = 1, 2$. There is a canonical $G$-equivariant quasi-conformal homeomorphism $\Omega : \Lambda_1 G \to \Lambda_2 G$ which induces a bijection $\omega$ from $E_G X_1$ to $E_G X_2$. As follows from the above construction of Gromov, the bijection $\omega$ can be realized by a homeomorphism $E_G X_1 \to E_G X_2$ sending $\mathbb{R}$-orbits to $\mathbb{R}$-orbits. M.Bourdon [28] proves that the bijection $\omega$ can be realized by a homeomorphism $E_G X_1 \to E_G X_2$ commuting with the flows if and only if the quasi-conformal homeomorphism $\Omega : \Lambda_1 G \to \Lambda_2 G$ is conformal. (Recall that, as in Theorem 3.5, a map $f : (A, d) \to (B, d')$ between two metric spaces is *conformal* if for any $a_0 \in A$, the limit of $\frac{d'(f(a), f(a_0))}{d(a, a_0)}$ as $a$ approaches $a_0$ exists, is finite and non-zero.)

Geodesic flow on hyperbolic groups is a useful tool which is not sufficiently well understood at the moment. Many ideas and techniques from the Riemannian case carry over provided the geodesic flow has sufficiently strong uniqueness properties (see for example the work of M.Gromov on bounded cohomologies of hyperbolic groups). For example, V.Kaimanovich [110] shows that under suitable uniqueness assumptions the geodesic flow on a hyperbolic group has very good measure-theoretic and ergodic properties (such as Hopf dichotomy for harmonic measures coming from Markov operators) which are similar to the classical Hopf-Tsuji-Sullivan Theorem [109, 173] for negatively curved Riemannian manifolds. In the $CAT(-1)$-case and negatively curved Riemannian manifold case there are deep connections between ergodic properties of the geodesic flow and various rigidity properties of the space (see for example [14, 15, 196, 195, 94]).

However, establishing the existence of good geodesic flows for arbitrary hyperbolic groups proved difficult for classes more general than groups acting on $CAT(-1)$-spaces.

## 14. Conformal Dimension

We have seen that the visual metric on the boundary of a hyperbolic space is not canonical, while the quasiconformal and quasi-Möbius structures are canonical. Moreover, the quasiconformal and quasi-Möbius structures are preserved by quasi-isometries of hyperbolic spaces. This led P.Pansu [150] to the following important concept, which by construction is a quasi-isometry invariant.

**Definition 14.1** (Conformal dimension). Let $(X, d)$ be a proper hyperbolic space. Let $d_a$ be a visual metric on $\partial X$. Define *conformal dimension* of $\partial X$

as:

$conf\,dim\,\partial X := \inf\{$ Hausdorff dimension of $d'\,|d'$ is a metric on $\partial X$ which is

quasiconformally equivalent to $d_a\}$.

In view of Theorem 3.2 the above definition does not depend on the choice of $d_a$ since any two visual metrics on the boundary are quasiconformally equivalent. Moreover, if two hyperbolic metric spaces $X$ and $Y$ are quasi-isometric, then $conf\,dim\,\partial X = conf\,dim\,\partial Y$ [99],[27]. Thus conformal dimension of the boundary is a quasi-isometry invariant of the space.

A simple but important observation states that (see for example [27]):

**Proposition 14.2.** *Let $(X, d)$ be a proper hyperbolic space. Then*

1.

$$dim\,\partial X \leq conf\,dim\,\partial X$$

*where $dim\,\partial X$ is the topological dimension of $\partial X$.*

2. *If $(Y, d')$ is another hyperbolic space which is quasi-isometric to $(X, d)$ then $conf\,dim\,\partial X = conf\,dim\,\partial Y$.*

3. *Let $(Y, d')$ be another hyperbolic space and let $h : X \to Y$ be a quasi-isometric embedding. Then $conf\,dim\,\partial X \leq conf\,dim\,\partial Y$.*

Conformal dimension of the boundary appears to be a very useful tool for distinguishing hyperbolic groups up-to quasi-isometry. However, the big drawback in the current understanding of this notion is that conformal dimension is extremely difficult to compute or even to estimate. Usually it is possible only for very symmetric spaces, such as trees or buildings. For example, for the classical hyperbolic space $\mathbb{H}_K^n$ (where $K$ is the field of either real or complex numbers or quaternions or Cayley numbers) P.Pansu [150] verified that

$$conf\,dim\,\partial\mathbb{H}_K^n = nk + k - 2 = dim\,\partial\mathbb{H}_K^n + k - 1,$$

where $k$ is the topological dimension of $K$. Thus in these cases the conformal dimension is equal to the Hausdorff dimension with respect to the standard metric. Other spaces for which the exact value of the conformal dimension is given are Fuchsian buildings [33],[31]. These buildings are $CAT(-1)$-spaces and the boundary of such a building is homeomorphic to the Menger curve. The conformal dimension is then expressed in terms of the growth rate of the volume of balls of a metric on the corresponding Weyl group.

There are some other important situations where conformal dimension can be estimated successfully. Namely, M.Gromov [96] suggested the idea of using the so-called *round trees* for this purpose. Suppose $T$ is a rooted simplicial tree drawn in the $xy$-plane so that the root is the origin and the whole tree $T$ is contained in the half-plane $y \geq 0$. A *round tree* $R(T)$ is obtained by rotating $T$ around the $y$-axis (sometimes the rotation is taken by the angle of $\pi$ rather than $2\pi$). Each infinite ray $r$ in $T$ starting at the root spans a topological two-plane in $R(T)$. This two-plane is given the metric of $\mathbb{H}^2$ in the Poincaré disk model. The resulting space $R(T)$ is hyperbolic and even $CAT(-1)$. It turns out that one can give an estimate from below on the conformal dimension of $\partial R(T)$ in terms of the length of edges of $T$ and the degrees of vertices of $T$. This provides an estimate from below on the conformal dimension of $\partial X$ for any $X$ containing a quasi-isometric copy of $R(T)$. This idea of M.Gromov [96] was made precise by M.Bourdon [27]. Note that the boundary of $R(T)$ contains a family of circles (or semi-circles if the rotation was by $\pi$) which correspond to the boundaries of $\mathbb{H}^2$ round leaves. This family turns out to have some good conformal properties and the estimates of conformal dimension rely on applying to this family the following useful result of P.Pansu [151]:

**Proposition 14.3.** *Let $X$ be a hyperbolic space. We assume that :*

1. *There is a family of curves $\mathcal{F}$ on $\partial X$ with a probability measure $\nu$ on $\mathcal{F}$.*
2. *There is a metric $d$ on $\partial X$ which is quasiconformally equivalent to any visual metric, Alhfors-regular, of Hausdorff dimension > 1.*
3. *There is a constant $C \geq 1$, such that for every ball $B$ of radius $r$ in $(\partial X, d)$, we have*

   $$C^{-1}r^{(H \dim d)-1} \leq \nu\{\gamma \in \mathcal{F} \, ; \gamma \cap B \neq \emptyset\} \leq Cr^{(H \dim d)-1},$$

   *(here $H \dim d$ stands for the Hausdorff dimension of $(\partial X, d)$).*

*Then $conf \dim \partial X = H \dim(\partial X, d)$.*

We recall that a metric $d$ on a space $E$ for Hausdorff dimension $D$ is said to be *Alhfors $D$-regular* (or just Alhfors-regular) if its $D$-Hausdorff measure $H$ satisfies the following property: there exists a constant $C \geq 1$ such that for every ball $B$ of $(E, d)$ of radius $r$, we have:

$$C^{-1}r^D \leq H(B) \leq Cr^D.$$

M.Bourdon successfully applied these ideas to estimate conformal dimensions of certain $CAT(-1)$ polygonal complexes and show that certain families of such complexes contain infinitely many quasi-isometry types [27], [30], [33].

J.Tyson [188] recently computed precise conformal dimensions for an interesting class of metric spaces:

**Theorem 14.4.** *Let $(X, d)$ be a proper Ahlfors $Q$-regular space ( where $Q > 1$) which supports an $(1, Q)$-Poincaré inequality. Then $conf \dim(X, d) = Q$.*

Poincaré inequality is an integral property of a metric space together with a Borel measure motivated by Sobolev spaces considerations (see [33], [105], [128], [102] for the definition and basic facts regarding Poincaré inequalities). The $(1, 1)$-Poincaré inequality is the strongest Poincaré inequality possible. Basic examples of metric examples supporting a $(1, 1)$-Poincaré inequality include Euclidean spaces, non-compact Riemannian manifolds of nonnegative Ricci curvature and Carnot groups. Other examples are given by the boundary of the hyperbolic building $(I_{p,q})$ studied by M.Bourdon [30].

If $G$ is a one-ended word-hyperbolic group such that $\partial G$ has local cut-points then $\partial G$ does not admit a Poincaré inequality. This follows from the results of B.Bowditch and J.Heinonen-P.Koskela [104]. There are no known examples of one-ended hyperbolic groups $G$ such that $\partial G$ has no local cut-points and $\partial G$ does not admit a Poincaré inequality.

Some other recent results related to conformal dimension were obtained in [129, 21, 22, 103] and other sources.

The relationship between the Pansu conformal dimension of the boundary of a hyperbolic group $G$ and the Gromov conformal dimension of $G$ (defined in Section 10) is not completely understood. Both of these dimensions are quasi-isometry invariants of the group which deserve substantial further investigation.

We conclude this section with the following important result of M.Gromov [96] connecting Pansu conformal dimension with $l_p$-cohomology (see [96] for the definitions):

**Theorem 14.5** (Non-vanishing theorem). *Let $G$ be a non-elementary word-hyperbolic group. Then $\overline{l_p H}^1(G) \neq 0$ for any $p > \max(1, conf\, dim\, \partial G)$.*

## 15. Sullivan-Patterson measures at the boundary and Mostow-type rigidity

15.1. **Sullivan-Patterson measures.** Given a convex-cocompact group of isometries $\Gamma$ of the hyperbolic space $\mathbb{H}^n$, S.J.Patterson [155] for $n = 2$ and D.Sullivan [171] for $n \geq 2$ introduced the concept of a conformal measure on the limit set $\Lambda\Gamma$, also called a *Sullivan-Patterson measure*. This measure turned out to be very useful in the theory of Kleinian groups and found many ergodic and geometric applications (see [172], [156] for more information).

M.Coornaert [64], [65] generalized many of these results to the setting of word-hyperbolic groups and Gromov-hyperbolic spaces. Before proceeding with the definitions of Sullivan-Patterson measures, let us recall some basic notions of measure theory.

**Convention 15.1.** Let $(X, d)$ be a metric space.

1. A Borel measure $\mu$ on $X$ is called *regular* if for every measurable set $A$
$$\mu(A) = \inf\{\mu(B)|A \subseteq B, \text{ where } B \text{ is a Borel subset of } X\}.$$

2. A measure $\mu'$ is absolutely continuous with respect to the measure $\mu$ if $\mu(A) = 0$ implies $\mu'(A) = 0$ for any $\mu$-measurable set $A$.

3. Suppose $\mu'$ is absolutely continuous with respect to $\mu$. The *Radon-Nikodym derivative* $d\mu'/d\mu$ is a $\mu$-integrable function $g \in L^1(\mu)$ such that for any $f \in L^1(\mu')$ we have
$$\int_X f\, d\mu' = \int_X fg\, d\mu = \int_X f\frac{d\mu'}{d\mu}\, d\mu.$$

From now on we will assume that $(X, d)$ is a proper hyperbolic space. All measures on $\partial X$ will be assumed to be Borel and regular. Also, if $\Gamma$ is a group acting on $X$ by isometries, $\mu$ is a measure on $\partial X$ and $\gamma \in \Gamma$, we will define a measure $\gamma^* \mu$ by setting $\gamma^* \mu(A) := \mu(\gamma(A))$ for $A \subseteq \partial X$. We need to recall the following important notion:

**Definition 15.2** (Busemann function). Let $x_0 \in X$ be a base-point and let $r : [0, \infty) \to X$ be a geodesic ray. The *Busemann function* or *horofunction associated to $r$* is the function $h : X \to R$ defined as:

$$h(x) := \lim_{t \to \infty} (d(x, r(t)) - t) \text{ for every } x \in X.$$

For every $p \in \partial X$ choose a geodesic ray $r_p : [0, \infty) \to X$ representing $p$ and denote by $h_p$ the Busemann function of $r_p$. Let $a > 1$ be a visual parameter. For any $p \in \partial X$ and any isometry $\gamma$ of $X$ we put

$$j_\gamma(p) := a^{\Delta(p)} \text{ where } \Delta(p) = h_p(x_0) - h_p(\gamma^{-1} x_0).$$

Now we can introduce a notion of a quasi-conformal measure on $X$ which is in a certain sense similar to the notion of a quasi-conformal metric discussed earlier.

**Definition 15.3** (Sullivan-Patterson measure). Let $D \geq 0$ and let $\mu$ be a nonzero measure on $\partial X$ with finite total mass. The measure $\mu$ is said to be $\Gamma$-*quasi-conformal of dimension $D$* if the measures $\gamma^* \mu$ ( where $\gamma \in \Gamma$) are absolutely continuous with respect to one another and if there exists $C \geq 1$ such that for every $\gamma \in \Gamma$,

$$\frac{1}{C} j_\gamma^D \leq d(\gamma^* \mu)/d\mu \leq C j_\gamma^D \qquad \mu\text{-almost-everywhere}.$$

It can be shown that the above notion does not depend on the choice of $x_0$. Moreover, quasi-conformality of a measure does not depend on the visual parameter $a$, although the dimension $D$ does. The main examples of quasi-conformal measures are Hausdorff measures. We recall the definition of a Hausdorff measure on a metric space.

**Definition 15.4.** Let $(E, d)$ be a metric space and $D \geq 0$. We denote $|U|$ the diameter of $U \subseteq E$. Given $\varepsilon \geq 0$, a family $(U_i)$ of subsets of $E$ is an $\varepsilon$-covering of $U$ if $U \subset \cup_i U_i$ and $|U_i| \leq \varepsilon$ for all $i$. We define for $U \subseteq E$

$$H_\varepsilon^D(U) = \inf \{ \underset{i}{\Sigma} |U_i|^D : (U_i) \text{ is an } \varepsilon\text{-covering of } U \}$$

and

$$H^D(U) = \lim_{\varepsilon \to 0} H_\varepsilon^D(U).$$

The number $H^D(U)$ is called the *$D$-dimensional Hausdorff measure* of $U$.

It turns out that for each $U$ there is a unique (possibly infinite) number $D_0$ such that $H^D(U) = \infty$ for any $D < D_0$ and $H^D(U) = 0$ for any $D > D_0$. This number $D_0$ is called the *Hausdorff dimension of $U$*.

An important observation of M.Coornaert, generalizing the results of D.Sullivan is:

**Proposition 15.5.** [65] *Let* $(X, d)$ *be a proper hyperbolic space and let* $d_a$ *be a visual metric on* $\partial X$ *with respect to a visual parameter* $a > 1$. *Let* $\Gamma$ *be a group acting by isometries on* $X$ *and let* $E \subset \partial X$ *be a* $\Gamma$*-invariant Borel subset (e.g. the limit set of* $\Gamma$*) with a finite nonzero Hausdorff* $D$*-measure with respect to the metric* $d_a$ *on* $\partial X$ *(thus* $D$ *must be the Hausdorff dimension of* $E$ *with respect to* $d_a$*). Then the Hausdorff* $D$*-measure* $H^D$ *on* $E$ *defines a* $\Gamma$*-quasi-conformal measure of dimension* $D$ *on* $\partial X$ *with respect to the visual parameter* $a$.

D.Sullivan originally used Hausdorff measures on limit sets of Kleinian groups to compute critical exponents of those groups. Critical exponent is an important geometric invariant which, in a sense, measures the degree of exponential growth of an orbit of the group.

**Definition 15.6** (Critical exponent). Let $\Gamma$ be a discrete group acting properly discontinuously by isometries on a hyperbolic space $(X, d)$. Let $x \in X$ be a base-point. The *critical exponent* $\delta_a(\Gamma)$ *with respect to a visual parameter* $a$ is defined as

$$\delta_a(\Gamma) := \lim_{R \to \infty} \sup \frac{log_a n(R)}{R}$$

where $n(R)$ is the number of points in the orbit $\Gamma x$ at distance at most $R$ from $x$. The critical exponent can be shown to be independent of the choice of $x$.

The critical exponent can also be characterized in terms of the Poincaré series. Namely, $\delta_a(\Gamma)$ has the property that the series

$$g_s(x) = \sum_{\gamma \in \Gamma} a^{-sd(x, \gamma x)}$$

diverges for $s < \delta_a(\Gamma)$ and converges for $s > \delta_a(\Gamma)$ (where $x \in X$ is a base-point). Accordingly, the group $\Gamma$ is said to be of $a$-convergent or $a$-divergent type depending on whether the Poincar'e series converges or diverges for $s = \delta_a(\Gamma)$. If $\Gamma$ is of divergent type and $\delta_a(\Gamma)$ is finite and nonzero, $\Gamma$ is called a *divergence group* with respect to the visual parameter $a$. For classical hyperbolic spaces $\mathbb{H}^n$ and $CAT(-1)$-spaces one almost always uses the visual parameter $a = e$ and $\delta(\Gamma)$ stands for $\delta_e(\Gamma)$.

Recall that an infinite group $\Gamma$ acting by isometries on a proper hyperbolic space $(X, d)$ is called *quasiconvex-cocompact* if the action is properly discontinuous, $\Gamma$ does not fix a point in $\partial X$ and for some $\Gamma$-invariant quasiconvex subset $A \subseteq X$ the quotient $A/\Gamma$ is compact. One can characterize convex-cocompact groups as follows:

**Lemma 15.7.** *Let* $\Gamma$ *be an infinite group acting properly discontinuously by isometries on a proper hyperbolic space* $(X, d)$ *so that there is no point in* $\partial X$ *fixed by* $\Gamma$ *(this implies that the limit set of* $\Gamma$ *contains at least two points). Then the following are equivalent:*

1. $\Gamma$ *is quasiconvex-cocompact.*

2. *Let $C(\Gamma)$ be the convex hull of the limit set $\Lambda\Gamma$ of $\Gamma$ (Recall that this means that $C(\Gamma)$ is the union of all biinfinite geodesics in $X$ with both endpoints in $\Lambda\Gamma$). Then $C(\Gamma)/\Gamma$ has finite diameter.*

*Proof.* We will sketch the idea of the proof. If (2) holds then for some bounded set $B \subseteq C(\Gamma)$ we have $\Gamma B = C(\Gamma)$. Since $X$ is proper, the closure $\overline{B}$ is compact. Since $C(\Gamma)$ is quasiconvex, it is easy to show that the set $A = \Gamma\overline{B}$ is quasiconvex and $\Gamma$-invariant. The quotient $A/\Gamma$ is compact and so $\Gamma$ is convex-cocompact.

Suppose now (1) holds and for some quasiconvex $\Gamma$-invariant subset $A \subseteq X$ the quotient $A/\Gamma$ is compact. Then the limit set $\Lambda A$ of $A$ in $\partial X$, which is defined as the collection of limits in $\partial X$ of sequences from $A$, can be shown to be closed. Since $\Lambda A$ is obviously $\Lambda$-invariant, this means that $\Lambda\Gamma \subseteq \Lambda A$. Since both $A$ and $C(\Gamma)$ are quasiconvex, we can show that $C(\Gamma)$ is contained in a bounded neighborhood of $A$. Since $A/\Gamma$ is compact, this implies that $C(\Gamma)/\Gamma$ is bounded, as required. $\square$

We have seen earlier that non-elementary quasiconvex subgroups of hyperbolic groups are quasiconvex-cocompact with respect to the action on the Cayley graph of an ambient hyperbolic group. Also, if $\Gamma$ is a discrete geometrically finite group of isometries of $\mathbb{H}^n$ without parabolics, then $\Gamma$ is quasiconvex-cocompact. Quasiconvex-cocompact groups are well-known to be of divergent type. The following theorem of M.Coornaert generalizes similar results of D.Sullivan for Kleinian groups.

**Theorem 15.8.** [65] *Let $(X, d)$ be a proper hyperbolic space and let $\Gamma$ be a group acting on $X$ by isometries properly discontinuously and quasiconvex-cocompactly. Let $d_a$ be a visual metric on $\partial X$ and let $\Lambda\Gamma$ be the limit set of $\Gamma$. Put $D = \delta_a(\Gamma)$ and let $H = H^D$ be the $D$-Hausdorff measure on $\Lambda\Gamma$ with respect to $d_a$.*
    *Then*

1. *The critical exponent $\delta_a(\Gamma)$ is equal to the Hausdorff dimension of $\Lambda\Gamma$ with respect to $d_a$. In particular $0 < H(\Lambda\Gamma) < \infty$.*
2. *The measure $H$ on $\Lambda\Gamma$ is $\Gamma$-quasiconformal of dimension $\delta_a(\Gamma)$.*
3. *The action of $\Gamma$ on the measure-space $(\Lambda\Gamma, H)$ is ergodic. This means that for any $\Gamma$-invariant Borel subset $B \subseteq \Lambda\Gamma$ either $H(B) = 0$ or $H(\Lambda\Gamma - B) = 0$.*
4. *If $\mu$ is a $\Gamma$-quasi-conformal measure of dimension $D'$ with support in $\Lambda\Gamma$, then $D' = D = \delta_a(\Gamma)$ and the measures $\mu$ and $H$ are absolutely continuous with respect to one another.*
5. *The measures $\psi \cdot H$, for any $\psi \in L_0^\infty(\Gamma)$, are all $\Gamma$-quasi-conformal with support $\Lambda\Gamma$ and they are the only $\Gamma$-quasi-conformal measures with support in $\Lambda\Gamma$. (Here $L_0^\infty(\Gamma)$ is the collection of all $H$-measurable functions on $\Lambda\Gamma$ which are pinched between two constant positive functions almost everywhere with respect to $H$.)*

It turns out that if $(X, d)$ is not just hyperbolic but a tree or a $CAT(-1)$-space, one can strengthen the above results and push further the analogy with Kleinian groups [68], [69], [29], [27], [47], [70]. This often involves ergodic-theoretic considerations similar to the work of D.Sullivan on the spectral theory of the Laplacian on the underlying hyperbolic space. Moreover, the restriction that $\Gamma$ be quasiconvex-cocompact can be relaxed, for example by considering measures on the conical limit set of $\Gamma$ instead of the full limit set [157], [158]. If $\Gamma$ is not convex-cocompact, Sullivan-Patterson measures do not have to be Hausdorff measures, but a close relationship between them remains (see for example [156], [69], [70]). F.Paulin has offered a version of Sullivan-Patterson theory for groups acting on hyperbolic spaces which is slightly different from the one described above [157], [158]. An inependent and detailed treatment of Sullivan-Patterson measures in $CAT(-1)$ case was also given by M.Burger and S.Mozes [47].

One can also approach Sullivan-Patterson measures from the point of view of symbolic dynamics. Namely, Sullivan-Patterson measures can often be represented as images of Gibbs measures on a certain shift space representing the symbolic dynamics on the limit set (the reader is referred to [41] for the background information on Gibbs measures). This was originally done by C.Series [169] for fuchsian groups (see also [131] for the free group case). M.Bourdon [26] extended these results to $CAT(-1)$-spaces. Namely, given a discrete quasiconvex-cocompact group $\Gamma$ acting by isometries on a proper $CAT(-1)$-space $X$. In this case M.Bourdon constructs a "coding map" $q : \Sigma' \to \Lambda\Gamma$ where $\Sigma'$ is a subshift of finite type in a shift $\Sigma(\mathbb{N}, S)$, where $S$ is a certain finite set (see Section 16 for precise definitions of Bernoulli shifts and subshifts of finite type). This map $q$ is surjective, Holder continuous and has uniformly bounded finite fibers. Then an explicit Gibbs measure $\mu$ is constructed on $\Sigma'$. It turns out that Sullivan-Patterson measures on $\Lambda\Gamma$ are connected with the Gibbs measure $\mu$:

**Theorem 15.9.** *In the above notations let $\delta = \delta(\Gamma)$ be the critical exponent of $\Gamma$ (so that $\delta$ is equal to the Hausdorff dimension of $\Lambda\Gamma$). Let $H$ be the Hausdorff measure of dimension $\delta$ on $\Lambda\Gamma$. Then the measures $H$ and $q_*\mu$ are absolutely continuous with respect to each other and their mutual densities are bounded.*

15.2. **Mostow-type rigidity.** Recall that by a famous result of D.Mostow if $\Gamma_1$ and $\Gamma_2$ are uniform lattices in $SO(n, 1) = Isom(\mathbb{H}^n)$, where $n > 2$, and if $h : \Gamma_1 \to \Gamma_2$ is an isomorphism then $h$ extends to a conjugation in $SO(n, 1)$. Some of the most interesting applications of Sullivan-Patterson measures are various rigidity results, generalizing the classical Mostow rigidity theorem. Roughly, the idea is the following. Suppose $X$ and $Y$ are two "particularly good" hyperbolic spaces. Then under suitable assumptions one can show that every topological embedding $\hat{f} : \partial X \to \partial Y$ which preserves cross-ratio, extends to an isometric embedding (and not just a quasi-isometric embedding) $f : X \to Y$. It often turns out that if $\Gamma_1 \leq Isom(X)$ and $\Gamma_2 \leq Isom(Y)$ so that $\Lambda\Gamma_1 = \partial X$ and $h : \Gamma_1 \to \Gamma_2$ is an isomorphism, and if some additional assumptions regarding the Sullivan-Patterson measures of $\Gamma_1$ and $\Gamma_2$ are satisfied, then the extension of $h$ to the boundary $\hat{h} : \Lambda\Gamma_1 \to \Lambda\Gamma_2$ is Möbius, that is preserves cross-ratio. This in turn implies that an isomorphism $h : \Gamma_1 \to \Gamma_2$ induces an isometric embedding $X \to Y$ and therefore a map $Isom(X) \to Isom(Y)$ extending $h$. The literature on Mostow-type rigidity for symmetric spaces is vast and it is far beyond the scope of this paper to survey it (we refer the reader to [127], [196], [152], [100], [136], [198] for more information). Therefore we will only mention the results which apply to "most" Gromov-hyperbolic spaces, as opposed to specific examples of such spaces.

M.Bourdon [29] has applied this approach to maps between uniform lattices in rank-one symmetric spaces and a convex-cocompact group of isometries of a CAT(-1)-space. Later his results were generalized by S.Hersonsky and F.Paulin [106] who obtained the following result.

**Theorem 15.10.** [106] *Let $X_1$ be a rank-one symmetric space, different from $\mathbb{H}^2$, with metric normalized so that maximal sectional curvature is $-1$ or let $X_1$ be a hyperbolic Bruhat-Tits building. Let $X_2$ be a proper CAT(-1)-space. Let $\Gamma_1$ and $\Gamma_2$ be discrete groups of isometries of $X_1$ and $X_2$ with the same critical exponents (with respect to the visual parameter $a = e$). Suppose $\Lambda\Gamma_1 = \partial X_1$ and that $\Gamma_2$ is a divergence group. Let $\hat{h} : \partial X_1 \to \partial X_2$ be a Borel map which is non-singular with respect to the Patterson-Sullivan measures and equivariant with respect to some morphism $h : \Gamma_1 \to \Gamma_2$. Then $\hat{h} : \partial X_1 \to \partial X_2$ is induced by an isometric embedding $X_1 \to X_2$.*

The crucial step in the proof of the above theorem is to show that under the assumptions of the above theorem the map $\hat{h} : \partial X_1 \to \partial X_2$ is a Möbius embedding, that is it preserves cross-ratio. This follows from the following general statement:

**Theorem 15.11.** [106] *Let $X_1$ and $X_2$ be proper CAT(-1)-spaces. Let $\Gamma_1$ and $\Gamma_2$ be discrete groups of isometries of $X_1$ and $X_2$ with the same critical exponents (with respect to the visual parameter $a = e$). Suppose that $\Gamma_2$ is a divergence group. Let $\hat{h} : \partial X_1 \to \partial X_2$ be a Borel map which is non-singular with respect to the Patterson-Sullivan measures and equivariant with respect to some morphism $h : \Gamma_1 \to \Gamma_2$. Then $\hat{h} : \partial X_1 \to \partial X_2$ is Möbius on the limit set of $\Gamma_1$.*

Note that if $X_1$ is a rank-one symmetric space and $\Gamma_1$ is a uniform lattice in the isometry group of $X_1$ then $\Gamma_1$ acts geometrically on $X_1$ and so $\partial X_1 = \Lambda\Gamma_1$. Also, if $\Gamma_2$ is a convex-cocompact group acting by isometries on $X_2$ then $\Gamma_2$ is a divergence group. Moreover, under these assumptions on $\Gamma_1$ and $\Gamma_2$ the Patterson-Sullivan measures are in fact Hausdorff measures of dimensions equal to the critical exponents of $\Gamma_1$ and $\Gamma_2$. In [106] S.Hersonsky and F.Paulin give many examples of divergence groups which are not necessarily convex-cocompact.

M.Burger and S.Mozes [47] obtained a number of Mostow-type rigidity results and results similar to Margulis superrigidity for discrete groups of isometries of $CAT(-1)$-spaces. They gave a self-contained treatment of Sullivan-Patterson measures for $CAT(-1)$ spaces and showed the relevance of ergodic-theoretic methods of G.Margulis [136], [198] and Poisson boundary considerations in this case.

We give an example of the many strong results obtained in [47]:

**Definition 15.12.** Let $\Gamma$ be a locally compact group acting by Borel automorphisms on a Borel space $B$ preserving a $\sigma$-finite measure class $\mu$.

The pair $(B, \mu)$ is called a *weak $\Gamma$-boundary* if

1. For any continuous affine action of $\Gamma$ on a separable locally convex space $E$ and any $\Gamma$-invariant compact convex subset $A$ of $E$ there is a $\mu$-measurable $\Gamma$-equivariant map $B \to A$.
2. The diagonal action of $\Gamma$ on $(B \times B, \mu \times \mu)$ is ergodic.

If $\Gamma$ is a closed subgroup of a locally compact group $G$ then a weak $\Gamma$-boundary $(B, \mu)$ is said to be a *weak $(G, \Gamma)$-boundary* if the action of $\Gamma$ on $B$ extends to a measure class preserving Borel $G$-action.

**Theorem 15.13.** [47]

*Let $X$ be a proper $CAT(-1)$-space and let $\Gamma \le Isom(X)$ be a discrete divergence group. Denote by $G$ the closure in $Isom(X)$ of the commensurator of $\Gamma$.*

*Then:*

1. *The Patterson-Sullivan measure class $\mu$ of $\Gamma$ on $\partial X$ is $G$-invariant and the pair $(\partial X, \mu)$ is a weak $(G, \Gamma')$-boundary for any subgroup of finite index $\Gamma' \le \Gamma$.*
2. *Let $Y$ be a proper $CAT(-1)$-space and let $\pi : \Gamma \to Isom(Y)$ be a homomorphism with a non-elementary image. Then there exists a unique $\Gamma$-equivariant measurable map $\phi : \partial X \to \partial Y$ and almost all values of $\phi$ lie in the limit set $\Lambda\pi(\Gamma) \subseteq \partial Y$.*
3. *Suppose $Y$ is as in (2) and $H \le Isom(X)$ be such that $\Gamma \le H \le G$. Suppose also that $\pi : H \to Isom(Y)$ is a homomorphism with a non-elementary image and such that the action of $\pi(H)$ on $\partial Y$ is c-minimal (that is $Y$ is the only nonempty closed convex $\pi(H)$-invariant subset of $Y$). Then $\pi$ extends to a continuous homomorphism*

$$\pi : \overline{H} \to Isom(Y).$$

We have already mentioned in Section 8 rather remarkable results of M.Bonk and B.Kleiner [24] regarding group actions on Ahlfors regular spaces. In the same paper M.Bonk and B.Kleiner also obtain the following strong rigidity result:

**Theorem 15.14.** *Let $k > 1$ and suppose $G$ is a group acting geometrically on a $CAT(-1)$-space $X$. Suppose both the topological and the Hausdorff dimensions of $\partial X$ are equal to $k$ (hence the conformal dimension of $\partial X$ is* *equal to k as well). Then X contains a convex G-invariant subset isometric to $\mathbb{H}^{k+1}$. In particular G acts geometrically on $\mathbb{H}^{k+1}$ and the boundary of both X and G is the k-sphere.*

Some other interesting rigidity results for groups acting on $CAT(-1)$-spaces were obtained by Y.Gao [89]. Most results which we quoted in this section can be substantially strengthened for isometry groups of special classes of hyperbolic spaces, such as the classical hyperbolic spaces $\mathbb{H}^n$, trees etc. We refer the reader to the work of C.Bishop-P.Jones [20], C.Yue [196] and M.Bonk-B.Kleiner [24] for a sample of such statetments.

## 16. SYMBOLIC DYNAMICS FOR THE ACTION ON THE BOUNDARY

Since a word-hyperbolic group G acts on $\partial G$ by homeomorphisms, the pair $(\partial G, G)$ can be considered as an abstract dynamical system and studied by means of symbolic dynamics. This approach was pursued by M.Coornaert and A.Papadopoulos in an excellent book [67] following the ideas of M.Gromov [95]. The symbolic dynamics point of view for studying boundaries of various classes of hyperbolic groups also has been pursued by many other authors (see [67] for a more comprehensive reference list).

Before stating the main results of [67] let us recall some basic notions of symbolic dynamics. A *dynamical system* is a pair $(\Omega, \Gamma)$ where $\Omega$ is a topological space and $\Gamma$ is a countable semigroup with a continuous action on $\Omega$. Note that if $\Gamma$ is not a group, individual elements of $\Gamma$ need not act by homeomorphisms on $\Omega$. If $Y \subseteq \Omega$ and $\gamma \in \Gamma$, we will denote

$$\gamma^{-1}Y := \{x \in \Omega | \gamma x \in Y\}.$$

**Convention 16.1.** Let $\Gamma$ be a countable semigroup and let $S$ be a finite set of symbols. Denote by $\Sigma = \Sigma(\Gamma, S)$ the set of all functions $\sigma : \Gamma \to S$. Since $\Gamma$ is countable, the set $\Sigma(\Gamma, S)$ can be thought of as the collection of sequences of symbols from $S$ indexed by elements of $\Gamma$. The semigroup $\Gamma$ and the set $S$ are equipped with discrete topology and the set $\Sigma(\Gamma, S)$ is given the product topology. Thus if $S$ is finite with at least two elements and if $\Gamma$ is infinite, then $\Sigma(\Gamma, S)$ is homeomorphic to the Cantor set. The semigroup $\Gamma$ has a natural left continuous action on $\Sigma(\Gamma, S)$:

$$(\gamma \cdot \sigma)(\gamma') := \sigma(\gamma'\gamma), \text{ where } \gamma, \gamma' \in \Gamma, \sigma \in \Sigma(\Gamma, S).$$

The dynamical system $(\Sigma(\Gamma, S), \Gamma)$ is called the *Bernoulli shift associated to $\Gamma$ on the set of symbols S*. By abuse of notation we will also often refer to $\Sigma = \Sigma(\Gamma, S)$ as a Bernoulli shift.

Bernoulli shifts and their subshifts play the role of free groups in symbolic dynamics. They can be used to encode more complicated dynamical systems, which leads to the notions of dynamical systems of finite type (parallel to finitely generated groups) and of finitely presentable dynamical systems.

**Definition 16.2** (Subshift of finite type). Let $\Sigma(\Gamma, S)$ be a Bernoulli shift. A *subshift* of $\Sigma$ is a closed $\Gamma$-invariant subset of $\Sigma$.

A subset $C \subseteq \Sigma$ is called a *cylinder* if there is a finite set $F \subseteq \Gamma$ and a collection $A$ of maps from $F$ to $S$ such that

$$C = \{\sigma \in \Sigma \,|\, \sigma|_F \in A\}.$$

A subshift $\Phi$ of $\Sigma$ is said to be *subshift of finite type* if there is a cylinder $C$ such that

$$\Phi = \bigcap_{\gamma \in \Gamma} \gamma^{-1} C.$$

**Definition 16.3** (System of finite type). A dynamical system $(\Omega, \Gamma)$ is said to be a *system of finite type* if there exist a finite set $S$, a subshift of finite type $\Phi \subseteq \Sigma(\Gamma, S)$ and a continuous, surjective and $\Gamma$-equivariant map $\pi : \Phi \to \Omega$.

Suppose $(\Omega', \Gamma)$ and $(\Omega, \Gamma)$ are two dynamical systems and $\pi : \Omega' \to \Omega$ is a surjective $\Gamma$-equivariant map. We will denote

$$R(\pi) := \{(x, y) \in \Omega' \times \Omega' \,|\, \pi(x) = \pi(y)\}.$$

Thus $R(\pi)$ is the graph of the equivalence relation on $\Omega'$. It is easy to see that $R(\pi)$ is a closed $\Gamma$-invariant subset of $\Omega' \times \Omega'$ and that $\pi$ induces a $\Gamma$-equivariant homeomorphism between the quotient space $\Omega'/R(\pi)$ and $\Omega$. Note also that if $\Omega'$ is a subshift of $\Sigma(\Gamma, S)$ then $R(\pi)$ is a subshift of the Bernoulli shift $\Sigma(\Gamma, S) \times \Sigma(\Gamma, S) = \Sigma(\Gamma, S \times S)$.

**Definition 16.4** (Finitely presented dynamical system). A dynamical system $(\Omega, \Gamma)$ is said to be *finitely presentable* if there exist a finite set $S$, a subshift of finite type $\Phi \subseteq \Sigma(\Gamma, S)$ and a continuous, surjective and $\Gamma$-equivariant map $\pi : \Phi \to \Omega$ such that $R(\pi)$ is a subshift of finite type of the Bernoulli shift $\Sigma(\Gamma, S \times S)$.

We are now able to state the first main result of M.Coornaert and A.Papadopoulos [67] which follows the ideas of the book of M.Gromov [95]:

**Theorem 16.5.** *Let $G$ be a word-hyperbolic group (so that $G$ is also a semigroup). Then the dynamical system $(\partial G, G)$ is finitely presentable.*

M.Coornaert and A.Papadopoulos [67] give two different proofs of the above result. The first proof uses so-called *cocycles* on a hyperbolic space $(X, d)$, that is real-valued functions from a bounded neighborhood of the diagonal in $X \times X$ which satisfy some natural cocycle and quasiconvexity conditions. It turns out that any cocycle $\phi$ is "integrable" , that is of the form $\phi(x, y) = h(x) - h(y)$ for some function $h : X \to \mathbb{R}$. Moreover, the gradient lines of any cocycle $\phi$ converge to a unique point at infinity $\pi(\phi) \in \partial X$. The idea of coding the boundary of a hyperbolic space using cocycles was originally suggested by M.Gromov [95].

If $X$ is the Cayley graph of a word-hyperbolic group $G$, one then considers a variation of the cocycle construction and studies the set $\Phi_0$ of integer-valued cocycles on the vertex set of $X$ (that is on $G$). Similarly to the general case there is a map $\pi_0 : \Phi_0 \to \partial G$ corresponding to limits of gradient lines. Then M.Coornaert and A.Papadopoulos consider a ball $B$ of sufficiently big finite radius in $G$ centered around the identity and define the set of symbols $S$ to be the set of restrictions to $B \times B$ of cocycles from $\Phi_0$. Note that $S$ is finite. Let $\Sigma = \Sigma(G, S)$ be the corresponding Bernoulli shift. Then $\Phi_0$ naturally injects in $\Sigma$ via the map $P : \Phi_0 \to \Sigma$, $P(\phi)(g) := g\phi|_{B \times B}$. Thus $\Phi_0$ can be considered as a subset of $\Sigma$, which turns out to be a subshift of finite type. It is then possible to show that the map $\pi_0 : \Phi_0 \to \partial G$ provides a finite presentation for the dynamical system $(\partial G, G)$.

The second finite presentation for $(\partial G, G)$ is based on a different idea and uses as the set of symbols $S$ a collection of "relations" on $G$ which are also thought of as discrete quasigeodesic vector-fields.

M.Coornaert and A.Papadopoulos also show that the boundary of a torsion-free hyperbolic group is a *semi-Markovian* space. Semi-Markovian spaces are special types of compact metrizable spaces which in a certain sense can be constructed out of a finite amount of combinatorial data. Thus finite simplicial complexes and many "regular" fractal spaces (e.g. Cantor set, Sierpinski carpet, Menger curve etc) are semi-Markovian.

**Definition 16.6.** Let $S$ be a finite set and $\Sigma = \Sigma(\mathbb{N}, S)$ be the Bernoulli shift where $\mathbb{N}$ is considered as a semigroup with respect to addition.

A subset $\Psi \subseteq \Sigma$ is called *semi-Markovian* if there is a cylinder $C \subseteq \Sigma$ and a subshift of finite type $\Phi \subseteq \Sigma$ such that $\Psi = C \cap \Phi$.

A compact metrizable topological space $\Omega$ is called *semi-Markovian* if there exist a finite set $S$, a semi-Markovian subset $\Psi \subseteq \Sigma(\mathbb{N}, S)$ and a continuous surjective map $\pi : \Psi \to \Omega$ such that the associated equivalence relation

$$R(\pi) = \{(x, y) \in \Psi \times \Psi \mid \pi(x) = \pi(y)\}$$

is a semi-Markovian subset of $\Sigma(\mathbb{N}, S) \times \Sigma(\mathbb{N}, S) = \Sigma(\mathbb{N}, S \times S)$.

**Theorem 16.7.** [67] *Let $G$ be a torsion-free word-hyperbolic group. Then the boundary $\partial G$ is a semi-Markovian space.*

M.Gromov later observed that the "torsion-free" assumption in the above theorem can be dropped. The boundary of a hyperbolic group is rather explicitly encoded in the finite state automaton recognizing the set of lexicographically least geodesic words. Thus it is not very surprising that the boundary of a hyperbolic group turns out to be semi-Markovian. The precise argument involves analyzing the sequence of projection maps similar to that described in Theorem 2.26 and investigating the rooted simplicial locally finite tree associated to that sequence. It can be shown that a semi-Markovian space has finite topological dimension and that there are only countably many homeomorphism types of semi-Markovian spaces. This yields an alternative proof that the boundary of a hyperbolic group is finite-dimensional. Note that since word-hyperbolic groups are finitely presentable and there are only countably many finite group presentations, it is obvious that the number of homeomorphism types of boundaries of hyperbolic groups is countable. Since spheres of all dimensions occur as such boundaries, this number is in fact infinite countable.

## 17. MISCELLANEOUS

As we have seen in dimensions zero and one the topological type of boundaries of hyperbolic groups is rather restricted (the basic examples being the Cantor set, the circle, the Sierpinski carpet and the Menger curve). It is far less clear and much less well understood what happens in higher dimensions. We know, of course, that for any $n \geq 1$ the sphere $S^n$ can occur as the boundary of a hyperbolic groups. A.N.Dranishnikov [72] constructed examples of right-angled hyperbolic Coxeter groups having Pontryagin surfaces and two-dimensional Menger compacta as their boundaries. (Note that right-angled hyperbolic Coxeter groups are known to have virtual cohomological dimension at most four and hence their boundaries are at most three-dimensional). A two-dimensional universal Menger compactum is an analog of the Menger curve with topological dimension two. If $p$ is a prime number then the Pontryagin surface $\Pi_p$ is a two-dimensional compactum with some exotic cohomological properties. Namely, $dim_{\mathbb{Q}}\Pi_p = dim_{\mathbf{Z}_q}\Pi_p = 1$ and $dim_{\mathbf{Z}_p}\Pi_p = 2$, where $q$ is a prime and $q \neq p$. A.Dranishnikov [72] shows that for any prime $p$ the space $\Pi_p$ occurs as the boundary of a right-angled hyperbolic Coxeter groups. At the present it is not known if for $n \geq 3$ the $n$-dimensional universal Menger compactum $\mu^n$ can be the boundary of a hyperbolic group. Also classifying groups which have the two-sphere or the Sierpinski carpet as the boundary remains a challenging open problem.

One says that a finitely generated group $G$ is *quasi-isometrically rigid* if for any group $G'$ quasi-isometric to $G$ there is a finite normal subgroup $N \leq G'$ such that $G'/N$ is isomorphic to a subgroup of finite index in $G$. We have seen in Section 15 that many hyperbolic groups exhibit Mostow-type rigidity. Any two word-hyperbolic groups having the circle as the boundary are virtually fuchsian and thus abstractly commensurable (see Section 5 above). We have also observed earlier that any group quasi-isometric to a nonabelian free group $F$ of finite rank is itself virtually free and thus commensurable with $F$. Therefore hyperbolic surface groups are quasi-isometrically rigid.

However, most other well-understood examples of hyperbolic groups (unlike lattices in higher rank semi-simple Lie groups [78], [127], [79] [165]) do not possess strong quasi-isometric rigidity properties. For example, any two groups acting geometrically on $\mathbb{H}^3$ are quasi-isometric to each other but they do not have to be abstractly commensurable or even weakly commensurable. It seems plausible that "generic" (in the sense of [148], [58]) small cancellation groups with the Menger curve as the boundary are in fact quasi-isometrically rigid, but nothing is known on this question to date.

In [121] M.Kapovich and B.Kleiner constructed an interesting example of a *topologically rigid* hyperbolic group $G$. That is to say, any homeomorphism $\partial G \to \partial G$ is induced by the action of an element of $G$. The group $G$ in this example is non-elementary and has two-dimensional boundary. Any topologically rigid group $G$ has the following property. If $G_1$ is a group whose boundary is homeomorphic to $\partial G$ (e.g. if $G_1$ is quasi-isometric to $G$) then $G_1/N$ embeds into $G$ as a subgroup of finite index (where $N$ is the maximal finite normal subgroup of $G_1$).

There are many other worthwhile points of view on the boundary which we have not discussed. For example, we have mentioned that the language of geodesic words forms an automatic structure on a word-hyperbolic group. W.Neumann and M.Shapiro [146] described how to associate to any asynchronously automatic structure on a group a certain *boundary*. This boundary consists of equivalence classes of infinite words which can be "read" in the automaton recognizing the automatic language and carries a natural topology similar to the topology on the hyperbolic boundary. Not surprisingly, for a word-hyperbolic group the boundary of the automatic structure given by the language of geodesics coincides with the hyperbolic boundary of the group. Moreover, W.Neumann and M.Shapiro define a natural notion of equivalence among asynchronously automatic structures on a given group and show that the topological type of the boundary of a structure depends only on its equivalence class. Somewhat surprisingly, N.Brady [43] proved that if $G$ is a hyperbolic surface group then any asynchronously automatic structure on $G$ is equivalent to the geodesic automatic structure. In the original paper W.Neumann and M.Shapiro [146] showed that for the fundamental group of a closed hyperbolic manifold fibering over a circle the number of equivalence classes of the asynchronously automatic structures is infinite even modulo the action of the commensurator.

I.Kapovich and H.Short [117] showed that for a hyperbolic group $G$ the set of rational points in $\partial G$ can be treated as a rudimentary boundary and it already carries a great deal of information about the group.

Another interesting general result regarding the actions of hyperbolic groups on their boundaries is due to S.Adams [1]. He proved that the natural action of a word-hyperbolic group $G$ on the boundary $\partial G$ is amenable with respect to any quasi-invariant measure on $\partial G$.

Department of Mathematics, University of Illinois at Urbana-Champaign, 1409 West Green Street, Urbana, IL 61801, USA
*E-mail address*: `kapovich@math.uiuc.edu`

Dept. of Mathematics, New York City Technical College, 300 Jay Street Brooklyn, NY 11201, USA
*E-mail address*: `NBenakli@nyctc.cuny.edu`